\documentclass[a4paper]{amsart}
\usepackage[left=3cm, right=2.5cm, top=2.5cm, bottom=2.5cm]{geometry}   

\usepackage{dsfont}
\usepackage{amsmath}
\usepackage{mathtools}
\usepackage{amssymb,xcolor}
\usepackage{mathrsfs}
\usepackage[numbers]{natbib}
\usepackage[pdfencoding=auto]{hyperref}
\usepackage{subcaption}

\newtheorem{theorem}{Theorem}
\newtheorem{definition}{Definition}
\newtheorem{lemma}{Lemma}
\newtheorem{proposition}{Proposition}
\newtheorem{remark}{Remark}
\newtheorem{example}{Example}
\newtheorem{corollary}{Corollary}

\newcommand{\N}{\mathbb{N}}
\newcommand{\R}{\mathbb{R}}

\newcommand{\bx}{{\bf{x}}}

\newcommand{\bxe}[1]{ {\bf{x}}({#1}) }
\newcommand{\bX}{{\bf{X}}}

\newcommand{\bC}{{\bf{C}}}

\newcommand{\by}{{\bf{y}}}

\newcommand{\bU}{{\bf{U}}}
\newcommand{\bz}{{\bf{z}}}
\newcommand{\bu}{{\bf{u}}}

\newcommand{\bA}{{\bf{A}}}

\newcommand{\bn}{{\bf{n}}}
\newcommand{\bc}{{\bf{c}}}
\newcommand{\bg}{{\bf{g}}}
\newcommand{\e}[1]{ {\bf{e}}_{{#1}} }
\newcommand{\bnu}{{\boldsymbol{\nu}}}

\newcommand{\ind}[1]{{\mathds{1}_{\{#1\}} }} 
\newcommand{\dd}{\textnormal{d}} %
\renewcommand{\P}{\mathbb{P}}
\newcommand{\lbl}[3]{\rho_{{#2}}^{{#3}}\left({#1}\right)}  
\newcommand{\E}{\mathbb{E}}

\newcommand{\type}[2]{\mathsf{g}(#1,#2)} 
\newcommand{\dtype}{d} 

\definecolor{colorn}{RGB}{73, 156, 0}

\newcommand*{\affmark}[1][*]{\textsuperscript{#1}}

\title{Fixation and stationary times for the \texorpdfstring{$\Lambda$}{Lambda}-Wright-Fisher process}

\author{A. Blancas\affmark[1,A] \and A. Gonz\'alez Casanova\affmark[2,3,B] \and S. Hummel\affmark[4,C] \and S. Palau\affmark[5,D]}

\address{ \newline 
	\affmark[1]Department of Statistics, ITAM, Mexico City, Mexico \newline
	\affmark[2]School of Mathematical and Statistical Sciences, Arizona State University, Tempe, USA \newline 
	\affmark[3]Center for Mechanisms of Evolution, Biodesign Institute, Arizona State University, Tempe, USA \newline
        \affmark[4]Department of Health Science and Technology, ETH Zurich, Zurich, Switzerland \newline
	\affmark[5]Departamento de Probabilidad y Estad{\'i}stica, Instituto de Investigaci{\'o}n en Matem{\'a}ticas Aplicadas y en Sistemas, Universidad Nacional Aut{\'o}noma de M{\'e}xico, M{\'e}xico City, M{\'e}xico}
\email{\affmark[A]airam.blancas@itam.mx, \affmark[B]agonz591@asu.edu \affmark[C]shummel@hest.ethz.ch and \affmark[D]sandra@sigma.iimas.unam.mx}

\begin{document}

\begin{abstract}
We study the fixation and stationary behavior of the $\Lambda$-Wright--Fisher process with parent-independent mutation and finitely many types, a jump-diffusion model for allele frequency dynamics in large populations with potentially large offspring variance. Using a lookdown construction, we characterize the distribution of fixation times and the order of allele extinctions in the absence of mutations, and identify a strong stationary time in the presence of mutations. Our results include explicit expressions for the mean fixation and stationary times for the Wright--Fisher diffusion, and mean fixation times in the Beta-coalescent case. A key component of our approach is the analysis of the \emph{fixation line} introduced by Hénard in [Ann. Appl. Probab., 25:3007–3032, 2015]. 
We extend this process to incorporate mutation, providing a unified framework for studying both fixation and equilibrium behavior.
\end{abstract}

\maketitle


\noindent\textbf{MSC2020 subject classifications:}
92D15, 
60J90, 
60J95, 
60J25, 
60J27.

\section{Introduction}

A central problem in population genetics and evolutionary biology is to determine the time until a subset of alleles goes extinct. 
The first time at which only one allele type remains is called the \emph{fixation time}. 
If mutations occur, genetic variation may be maintained, and over time the allele frequencies may reach an equilibrium. 
In this case, we say that the system reaches \emph{stationarity}.

For very large, randomly reproducing populations comprising individuals that carry one of $(\dtype+1)$ allelic types, 
the \emph{$(\dtype+1)$-type $\Lambda$-Wright-Fisher process} with parent-independent mutation can serve as a formal mathematical modeling framework.
This is a jump-diffusion process on the $d$-dimensional simplex that is characterized by three parameters:
a measure~$\Lambda$ on $[0,1]$, 
which encodes the distribution of offspring sizes in reproduction events; 
a non-negative real number~$\theta$, which specifies the total mutation rate;
and a vector $\bnu$ in the $d$-dimensional simplex
\[
\Delta_d\coloneqq\left\{\bx=(\bxe{1},\bxe{2},\ldots,\bx(d))\in [0,1]^{d}: \sum_{i=1}^d  \bx(i)\leq 1 \right\},
\] 
where the $i$-th component of~$\bnu$ corresponds to the probability of a mutation resulting in type~$i$.
For $\bx\in \Delta_d$, we set $\bxe{\dtype+1}\coloneqq 1-\sum_{i=1}^d  \bxe{i}$.
We write $\mathrm{int}(\Delta_{\dtype})$ for the interior of $\Delta_{\dtype}$.

To formally describe the mathematical framework,
we require some notation. 
We set $\N = \{1,2,\dots\}$, 
$\N_0 = \N\cup\{0\}$, and $[d] \coloneqq \{1,\ldots,d\}$ for a fixed $d \in \N$.
Vectors and sequences are written in boldface. 
In particular, for $i \in [\dtype]$, we write $\e{i}$ for the $i$-th unit vector in $\R^d$, and by an abuse of notation $\e{\dtype+1} \coloneqq (0,\ldots,0)\in \R^d$.
The set of twice continuously differentiable functions on $\Delta_d$ is denoted by $\mathcal{C}^2(\Delta_d)$.

Fix $\dtype\in \N$, $\theta\geq 0$, 
$\bx, \bnu\in \Delta_{\dtype}$ and $\Lambda$ a finite Borel measure on $[0,1]$.
The $(\dtype+1)$-type \emph{$\Lambda$-Wright--Fisher process with parent-independent mutation} parametrised by $\theta$ and $\bnu$, is the Markov process $\bX=(\bX_t:\, t\geq 0)$ on~$\Delta_\dtype$ with infinitesimal generator given by 
\begin{equation}
\mathcal{L}f(\bx)=\mathcal{L}_{\mathrm{WF}}f(\bx)+\mathcal{L}_{\Lambda}f(\bx)+\mathcal{L}_{\theta}f(\bx), \qquad f\in \mathcal{C}^2(\Delta_\dtype),  \label{eq:wfgenerator}
\end{equation} where
\begin{align}
	\mathcal{L}_{\mathrm{WF}}f(\bx)\coloneqq&\frac{\Lambda(\{0\})}{2}\sum_{i,j=1}^\dtype \bxe{i}\Big(\ind{j}(i)-\bxe{j}\Big)\frac{\partial^2 f}{\partial \bxe{i}\partial \bxe{j}} ,\label{eq:wfterm}
    \end{align}
    
    \begin{align}
	\mathcal{L}_{\Lambda}f(\bx)\coloneqq&\int_{(0,1]}  \sum_{i=1}^{\dtype+1}\bx(i) \left[f\left((1-r)\bx+r\e{i}\right)-f\big(\bx\big)\right]\frac{\Lambda(\dd r)}{r^2},\label{eq:lambdaterm1} \\
	\mathcal{L}_{\theta}f(\bx)\coloneqq& \theta\sum_{i=1}^{\dtype}\left(\big(1-\bxe{i}\big)\bnu(i) - \bxe{i}\sum_{j=1,\, j\neq i}^{\dtype+1}\bnu(j)\right) \frac{\partial f}{\partial \bxe{i}}. \label{eq:mutterm}
\end{align}
We write  $\bX^{\bx}=(\bX^{\bx}_t:\, t\geq 0)$ when $\bX_0=\bx$ a.s.
The existence of a solution to the martingale problem associated with~$\mathcal{L}$ has been frequently discussed in the literature; see for example~\cite[Lem.~4.5]{Birkner2009}. 

The process $\bX^{\bx}$ describes the evolution of type frequencies in an infinite haploid population in which individuals carry one of the $(\dtype+1)$ genetic types.
For $i\in[\dtype+1]$ and $t\geq 0$, the value $\bX^{\bx}_t(i)$ is the frequency of type $i$ at time~$t$.
The term~\eqref{eq:wfterm} accounts for frequency fluctuations arising from (neutral) reproduction events with small offspring size, where $\Lambda({0})$ determines the occurence rate of these events.
The term \eqref{eq:lambdaterm1} models (neutral) offspring events in which the offspring is large and makes up a fraction~$r$ of the population, with $r$ chosen according to the measure $r^{-2}\Lambda(\dd r)$.
The type producing the offspring in small and large reproduction events is selected uniformly at random from the current population.
The term \eqref{eq:mutterm} models mutation, where each individual independently mutates at rate~$\theta$.
Upon mutation, the resulting type is $i$ with probability~$\bnu(i)$.
Thus, types different from~$i$ mutate to type~$i$ at rate $\theta\bnu(i)$, leading to an increase in the frequency of type~$i$. Whereas type $i$ mutates to any other type at rate $\theta\sum_{j=1}^{\dtype+1}\bnu(j)\ind{j\neq i}$, leading to a decrease in the frequency of type~$i$.

The Wright-Fisher diffusion corresponds to $\Lambda=\delta_0$ and $\theta=0$, and it is a classic and still common null model for population genetic data.
Explicit expressions for the mean fixation time were derived by Litter in ~\cite{Littler1975} 
(see also~\cite[Ch.~8]{Durrett2008} for more context).
Moreover, Littler derives an expression for the probability of extinction to occur in a given order.
His expressions are relevant for biological applications, for example, in estimating divergence between (human) populations from linkage disequilibrium~\cite{sved2008divergence} and for excluding a neutral null model in favor of balancing selection at a major histocompatibility complex locus in microtine rodents~\cite{Oliver2012}.

There are biological settings, when the Wright-Fisher diffusion approximation is an inappropriate model because of large offspring variance, such as for certain marine species~\cite{Eldon2006}.
Then general measures~$\Lambda$ are relevant.
If $\Lambda$ is a general measure and $\theta=0$, 
 results that quantify the fixation time exist for the case of two types  (e.g., \cite[Prop.~2.29]{Cordero2022}), albeit these results are not explicit. Moreover, we are not aware of previous results that determine extinction to occur in a given order.
When $\theta > 0$, all types can coexist, and one can instead study the time to reach stationarity. Also here, the few explicit results deal with the two-type case, either in the finite population~\cite{donnelly2000convergence} or classic Wright-Fisher diffusion setting~\cite{Tavare1984}.

In this manuscript, we formulate a probabilistic framework to address questions related to fixation times, extinction order and time to stationarity, 
based on the lookdown model introduced by Donnelly and Kurtz ~\cite{Donnelly1999} and the fixation lines defined in \cite{henard2015}. We first introduce the corresponding mathematical notions to quantify the time to fixation and the time to stationarity. For questions concerning fixation times, we denote by $T_{\mathrm{lost},i}^{\bx}$ the first time at which type $i \in [\dtype+1]$ disappears, and by $T_{\mathrm{fix},n}^{\bx}$ the first time at which no more than $n \in [\dtype+1]$ types remain alive, i.e.
\begin{align*}
	T_{\mathrm{lost},i}^{\bx}\coloneqq\inf\{t\geq 0: \bX_t^{\bx}(i)=0\}
\qquad \text{and}\qquad T_{\mathrm{fix},n}^{\bx}\coloneqq \min_{\substack{J\subset[\dtype+1]\\ \lvert J\vert\geq \dtype+1-n}} \max_{i\in J}T_{\mathrm{lost},i}^{\bx}.
\end{align*}
The concept of time to stationarity can be formalized with the notion of a \emph{strong stationary time}.
This special stopping time was introduced for Markov chains in~\cite{Diaconis1990} and extended to more general Markov processes~\cite{Fill2016,Miclo2017}.

The following definition is adapted from Miclo~\cite{Miclo2017}.
Consider an ergodic Markov process $X = (X_t:\, t \geq 0)$ defined on a filtered probability space $(\Omega, \mathcal{G}, P, (\mathcal{G}_t:\, t \geq 0))$ with invariant distribution~$\pi$.
A $\mathcal{G}$-stopping time $T$ taking values in $[0, \infty)$ is said to be strong (for~$X$) if $T$ and $X_T$ are independent.
It is called a strong stationary time (for~$X$) if, in addition, $X_T$ is distributed according to~$\pi$.

Next, the lookdown model is constructed in a probability space that carries a countable infinite number of independent homogeneous Poisson processes and a sequence $\bU=(\bU(i): i\in \N)$ of independent uniform random variables on $[0,1]$.
The uniform variables are identifiers for the individuals in the population alive at time~$0$.
Specifically, each individual is coded by an element in $[0,1] \times \Delta_{\dtype}$, with the first component being an identifier (ultimately carrying information about the individual's ancestor) and the second component being used to determine the individual's type via the \emph{(genetic) type function}~$\mathsf{g}$, which is defined as
\begin{equation}
	\mathsf{g}:[0,1]\times \Delta_{d}\to[\dtype+1],\ \quad \ (u,\bx)\mapsto 	\mathsf{g}(u,\bx):=\min\left\{i\in[\dtype+1]: \sum_{j=1}^i \bx(j)>u \right\}\label{eq:typingfct}.
\end{equation}
For $i \in [\dtype+1]$, define the index of the first uniform that gets assigned genetic type~$i$ under~$\bx$,
\begin{equation}
	\mathfrak{m}^{\bx}_i\coloneqq \min\{k:\type{\bU(k)}{\bx}=i\} \label{eq:defmx},
\end{equation}
with the convention that $\min \emptyset=\infty$.
In particular, if $\bx(i)=0$, then $\mathfrak{m}^{\bx}_i=\infty$.
Observe that for $i\in [\dtype+1]$ and $\bx$ such that $\bx(i)\in (0,1]$, $\mathfrak{m}^{\bx}_i$ is a geometric random variable with parameter $\bx(i)$.
Moreover, for $n\in[\dtype+1]$ we set  $V_n^{\bx}$ as the first index where the uniforms have assigned $n$ different types under ${\bx}$,
\begin{equation}\label{eq:defVktyp}
V_n^{\bx}:= \min\left\{\ell: \lvert\{\type{\bU(k)}{\bx}: k\in[\ell] \right\}\rvert=n\},
\end{equation}
with the convention that $\min\emptyset = \infty$. In particular, if $\lvert \{i:\bx(i)>0\}\rvert = m<n$, then $V_n^\bx =\infty$.

The lookdown model provides a convenient framework to track genealogical relationships between individuals across time and it has played a key role in recent developments in population genetics.
Within this model, the backward genealogy of a population can be described by $\Lambda$-coalescents, which are certain type of exchangeable coalescent processes that were independently introduced by Pitman~\cite{Pitman1999} and Sagitov~\cite{Sa99}. In the $\Lambda$-coalescents, given $n$ genealogical lines, any subset of $k$ lines coalesces independently into one at rate $ \lambda_{n,k}+\ind{k=2}\Lambda(\{0\})$, where
 \begin{equation}\label{eq:ratetransLam}
 \lambda_{n,k}\coloneqq \int_{(0,1]}r^{k-2}(1-r)^{n-k}\Lambda(\dd r), \qquad  2\leq k\leq n.     
 \end{equation}
For most of our results, we impose a condition on the measure~$\Lambda$: namely that the associated $\Lambda$-coalescent comes down from infinity, i.e. if we start with infinitely many genealogical lines, after any positive amount of time the number of ancestral lineages becomes finite.
A necessary and sufficient condition for this to occur was established by Schweinsberg~\cite{schweinsberg2000}, and it reads 
\begin{equation}\label{cond:cdi}
\sum_{n=2}^\infty \left( \sum_{k=2}^n (k-1) \binom{n}{k} [\lambda_{n,k}+\ind{k=2}\Lambda(\{0\})] \right)^{-1} < \infty. 
\end{equation}
In the present manuscript, our two main examples are the Wright-Fisher diffusion and the Beta coalescent, where $\Lambda$ follows a $\mathrm{Beta}$ distribution. 
The biological relevance of the Wright-Fisher diffusion has already been discussed above.
The $\mathrm{Beta}$ distribution setting carries biological significance, for example because it naturally arises in systems with a certain dormancy component, see~\cite{Cordero2022a} for details.
\begin{example}[Coming down from infinity]$ \ $ \label{example}
	\begin{enumerate}
		\item [(a)] The Wright-Fisher diffusion arises when $\Lambda=\delta_{\{0\}}\Lambda(\{0\})$. 
        The corresponding coalescent is the Kingman's coalescent and it directly follows from~\eqref{cond:cdi} that this coalescent comes down from infinity.
        Moreover, whenever there is a Kingman component, i.e. $\Lambda(\{0\})>0$, the associated $\Lambda$-coalescent always comes down from infinity.
		\item  [(b)] If $\Lambda$ has  $\mathrm{Beta}\ (2-\alpha,\alpha)$-distribution with $\alpha\in (1,2)$. i.e. it is given by  
\begin{equation}
\Lambda(\dd x)=\mathrm{Beta}(2-\alpha,\alpha)(\dd x) := \frac{x^{1-\alpha}(1-x)^{\alpha-1}}{\Gamma(2-\alpha)\Gamma(\alpha)}\ind{[0,1]}(x)\dd x, \label{eq:measurelambdabeta}
\end{equation} 
then the associated coalescent comes down from infinity, see e.g. Berestycki \cite[Cor.~3.2]{Berestycki2009}.
	\end{enumerate}
\end{example}

Finally, as an extension of previous works \cite{henard2015, Pfaffelhuber2006} we will define a notion of  fixation lines which is embedded into the previous setup and it is compatible with mutation. For $k\in \N_0$, the $k$-\emph{fixation line} $F^k=(F_t^k:\, t\geq 0)$ will be a $\N_0$-valued, non-decreasing Markov chain starting at $k$ that jumps from $n$ to $n+\ell$ at rate
\[\ind{\ell=1}\bigg[\Lambda(\{0\})\binom{n+1}{2}+\theta (n+1)\bigg]+\binom{n+\ell}{\ell+1}\int_{(0,1]} r^{\ell+1}(1-r)^{n}\frac{\Lambda(\dd r)}{r^2},\qquad \ell\geq 1.\]
Let $T_{\infty}(k)\coloneqq \inf\{t\geq 0: F_t^k=\infty\}$ be its {\it{explosion time}}, with the convention that $\inf\emptyset = \infty$.
We will notice that under Condition~\eqref{cond:cdi} for $k>0$, $T_{\infty}(k)<\infty$. Besides, if $\theta = 0$, then $T_{\infty}(0)=\infty$. We allow that the fixation line starts at $\infty$, in which case it remains at $\infty$, that is, $(F_{t}^{\infty})_{t\geq 0}\equiv\infty$ and $T_{\infty}(\infty)=0$.

Now we are in position to state our first result.

\begin{theorem}[Fixation and stationary time] \label{thm:fixandstattimes}
Assume Condition~\eqref{cond:cdi} holds. Let $n\in [\dtype+1]$ and $\bx\in \Delta_{\dtype}$.
    \begin{enumerate}
        \item[(i)] Let $\theta=0$. Then, \begin{enumerate}
            \item The first time the $n$-th type vanishes from the population has the same distribution 
            as the explosion time of the $(m-1)$-th fixation line, where  $m$  is the index of the first uniform that gets assigned genetic type~$n$, i.e.
            \[T_{\mathrm{lost},n}^{\bx}\overset{(d)}{=}T_{\infty}(\mathfrak{m}^{\bx}_n-1).\]
            \item The first time there are not more than $n$ types has the same distribution  as the explosion time of the $(k-1)$-th fixation line, where  $k$  is the first index of the uniforms where there are $n+1$ different  assigned types, i.e.
            \[T_{\mathrm{fix},n}^{\bx}\overset{(d)}{=} T_{\infty}(V_{n+1}^{\bx}-1).\]

        \end{enumerate} 
        \item[(ii)] Let $\theta>0$ and let $\bnu\in \mathrm{int}(\Delta_{\dtype})$. Then $\bX^{\bx}$ has a stationary time~$T$ and \[T\overset{(d)}{=} T_{\infty}(0).\]
        Moreover, if $\Lambda=\Lambda(\{0\})\delta_0$, then $T$ is a strong stationary time.
    \end{enumerate}
\end{theorem}
In (ii), if $\bnu$ is not in the interior of~$\Delta_{\dtype}$, some types may go extinct and mutations will not be able to maintain the full genetic diversity. 
The proof of part (i) and (ii) of the Theorem~\ref{thm:fixandstattimes} is provided in Sections~\ref{sec: fixtime} and~\ref{sec:proofmixingtimes}, respectively. Both proofs use the lookdown construction.

\begin{figure}[t]
	\begin{center}
		\scalebox{0.97}{\includegraphics[width=14cm]{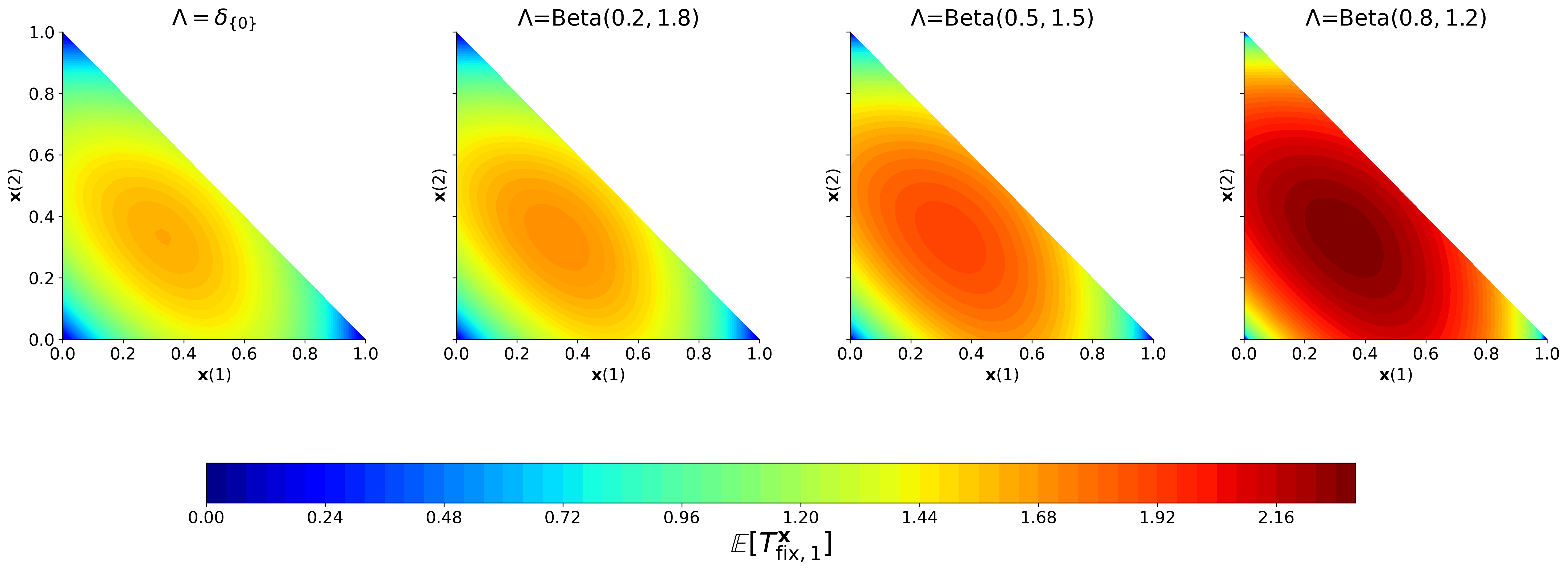}	}
	\end{center}
	\caption{Mean fixation time for three different types of alleles with initial frequency $\bx(1)$, $\bx(2)$, and $\bx(3)=1-\bx(2)-\bx(3)$. From left to right, $\Lambda$ is $\delta_0$ and $\mathrm{Beta}(2-\alpha,\alpha)$ for $\alpha=1.8$, $\alpha=1.5$, and $\alpha=1.2$. As the parameter~$\alpha$ increases, the mean fixation time increases for any~$\bx\in \Delta_{\dtype}$.}
	\label{fig:meanfixationtime}
\end{figure}

In the biologically relevant cases that are discussed in Example~\ref{example}, we can obtain more explicit results regarding the law of the explosion time. This leads to expressions for the mean fixation times for these examples. See Figure~\ref{fig:meanfixationtime} for an illustration of the mean fixation time with $n=1$ and three types of alleles. For the Wright--Fisher diffusion, we also provide the characteristic function of the fixation time, and in the subsequent result the mean strong stationary time.

\begin{corollary}[Mean fixation time]\label{coro:fixandstattimes_explicit}
Assume Condition~\eqref{cond:cdi} holds and that $\theta=0$. 
Let $n\in [\dtype]$ and $\bx\in \mathrm{int}(\Delta_{\dtype})$.
	\label{thm:meanfixk}$\ $
	\begin{enumerate}
		\item[(a)] \textbf{Wright-Fisher diffusion.} If $\Lambda=\Lambda(\{0\})\delta_0$, then
		\begin{align*}
		&	\E\left[T_{\mathrm{fix},n}^{\bx}\right]=\\
&\frac{-2}{\Lambda(\{0\})}\sum_{\ell=1}^n (-1)^{n-\ell}\binom{\dtype-\ell}{d-n}\sum_{1\leq c_1<\cdots<c_{\ell}\leq \dtype+1} \left(1-\sum_{j=1}^\ell \bx(c_j)\right)\log\left(1-\sum_{j=1}^\ell \bx(c_j)\right).
		\end{align*}
Moreover, then for any $t\in \R$,
	\begin{align*}
	 \E&\left[e^{it T_{\mathrm{fix},n}^{\bx}}\right]  =\sum_{\ell=1}^{n} (-1)^{n-\ell} \binom{d-\ell}{d-n} \\
    & \cdot \sum_{1\leq c_1<\ldots<c_\ell\leq d+1} \Bigg(1-{ \sum_{j=1}^{\ell}}\bx(c_j)\Bigg)\sum_{k=1}^{\infty}\Bigg(\sum_{j=1}^{\ell}\bx(c_j)\Bigg)^k\prod_{r=k}^\infty \Big(1-\frac{2it}{\Lambda(\{0\})(r+1)r} \Big)^{-1}.	
	\end{align*}
		\item[(b)] \textbf{Beta-coalescent.} If $\Lambda$ has  $\mathrm{Beta}\ (2-\alpha,\alpha)$-distribution with $\alpha\in (1,2)$, then
	\end{enumerate}
 \begin{equation*}\label{eq:firsttimekbeta}
		\begin{split}
			\E\left[T_{\mathrm{fix},n}^{\bx}\right]=\alpha&(\alpha-1)\sum_{\ell=1}^n(-1)^{n-\ell}\left\{\binom{\dtype-\ell}{d-n}\sum_{1\leq c_1<\cdots<c_{\ell}\leq\dtype+1} \left(1-\sum_{j=1}^\ell \bx(c_j)\right)\left(\sum_{j=1}^\ell \bx(c_j)\right)\right.  \\
			& \cdot\left.\int_{(0,1)} \frac{y(1-y)^{-1} }{\big(1-y\sum_{j=1}^\ell \bx(c_j) \big)\big((1-y)^{1-\alpha}-1\big)}\dd y\right\}.
			\end{split}
	\end{equation*}

\end{corollary}
 We remark that $\mathrm{int}(\Delta_{\dtype})$ is not a real restriction, as the Wright-Fisher process without mutation can always be reduced to that setting.
 The well-known mean fixation time formula of the Wright--Fisher diffusion is contained in our result, by taking $\Lambda=\delta_0$ and $n=1$, i.e.
\begin{equation*}
\E\left[T_{\mathrm{fix},1}^{\bx}\right]=-2\sum_{j=1}^{\dtype+1}\left(1-\bx(j)\right)\log\Big(1-\bx(j)\Big).
\end{equation*}
Littler \cite{Littler1975} (see also~\cite[Ch.~8.1.1]{Durrett2008}) derived the previous equality by using methods from diffusion theory that do not seem to translate easily to the non-diffusive case. In contrast, our proof uses probabilistic arguments which are detailed in Section~\ref{sec: fixtime}.

\begin{corollary}[Mean strong stationary time]\label{coro:meanstrongstationarytime}
    Assume Condition~\eqref{cond:cdi} holds, $\theta>0$ and $\bnu$ lies in the interior of~$\Delta_{\dtype}$.
    Moreoever, assume $\Lambda=\Lambda(\{0\})\delta_0$.
    Then $\bX^{\bx}$ has a strong stationary time~$T$ and
	\[
    \E\left[T\right]=\frac{2}{\Lambda(\{0\})} \sum_{j=1}^{\infty} \frac{1}{j\left(j-1+\frac{2\theta}{\Lambda(\{0\})}\right)}=
\frac{2}{\Lambda(\{0\})}\left(\psi\left(\frac{2\theta}{\Lambda(\{0\})}\right)+\gamma \right),
\] 
	where $\psi$ is the digamma function and $\gamma$ is the Euler-Mascheroni constant.
\end{corollary}
Without much additional effort, our lookdown construction also allows us to determine the probability for a given order of disappearance of types.
It turns out that this probability only depends on the initial frequencies, 
and is independent of any measure~$\Lambda$ that satisfies Condition~\eqref{cond:cdi}.
Denote by $\sigma_{[d+1]}$ the set of permutation of $[\dtype+1]$.

\begin{proposition}[Order of disappearance]\label{prop:disappearanceorder}
	Assume $\theta=0$ and Condition~\eqref{cond:cdi} holds. Let $(c_1,\ldots,c_{\dtype+1})$ be an element in $\sigma_{[d+1]}$. Then, 
    \[
\P(T_{\mathrm{lost},{c_{\dtype+1}}}^{\bx}<\cdots<T_{\mathrm{lost},{c_{1}}}^{\bx})
    =\prod_{i=1}^{\dtype+1} \frac{\bx(c_i)}{1-\sum_{j=1}^{i-1}\bx(c_j)}.
    \]
\end{proposition}

Proposition~\ref{prop:disappearanceorder} is proved in Section \ref{sec: fixtime}. By summing over all permutations such that $c_{d+1}=\eta\in [\dtype+1]$, we obtain the following corollary, which provides an explicit formula for the probability of the first type to disappear. This probability, again, is independent of~$\Lambda$ and just depends on the initial frequencies.  

\begin{corollary}[First type to disappear]
	Assume $\theta=0$ and Condition~\eqref{cond:cdi} holds. For $\eta\in [\dtype+1]$,
	\[\P( \text{Type $\eta$ is the first to disappear} )= \sum_{
	(c_1,\ldots, c_{d}, \eta)\in \sigma_{[d+1]}
	} \, \prod_{i=1}^{\dtype} \frac{\bx(c_i)}{1-\sum_{j=1}^{i-1} \bx(c_j)}.
\]
\end{corollary}
Note that by taking $d=3$ and $\eta=3$ in the corollary above, we obtain 
\[
\bx(1)\bx(2)((1-\bx(1))^{-1}+(1-\bx(2))^{-1}),
\]
which gives the probability that type~$3$ is the first to disappear. This is a well-known result in the case $\Lambda = \delta_0$ (see, e.g.,~\cite[Thm.~8.1]{Durrett2008}),
but has, to the best of our knowledge,
not been proved in more general cases.

The remainder of the manuscript is organized as follows.
Section~\ref{sec:pathwise_construction_lookdown} contains details on the lookdown model, 
its convergence to the $\Lambda$-Wright-Fisher process,
and a description of a dual process.
The fixation line process embedded into the lookdown construction is rigorously defined in Section~\ref{sec:fixation_line_process};
there we also define formally its explosion time and derive explicit expression for its mean in the case of the Wright-Fisher diffusion and the Beta-coalescent.
Finally,
our results about fixation and stationary times are proved in Section~\ref{sec:proof_fixation_and_stationary_times}.

\section{Pathwise construction and dual process}\label{sec:pathwise_construction_lookdown}
The process $\bX^\bx$ is formally defined as the solution to a martingale problem.
Although techniques from stochastic differential equations could potentially be used to study questions such as fixation time and time to stationarity,
we instead adopt a probabilistic framework based on the \emph{lookdown construction} of Donnelly and Kurtz~\cite{Donnelly1999}.
In this section, we develop a modification of the lookdown construction and show, using a dual process, that the empirical distributions of types converge to the process $\bX^\bx$.
The main difference to previous works is a novel implementation of mutation, which proves to be especially effective for investigating issues related to the distance to stationarity.

\subsection{Lookdown construction}\label{seq:lookdownconstruction}

The lookdown construction facilitates keeping track of the genealogy and the type of the individuals throughout time.
Individuals alive at time $t \geq 0$ are positioned at $(t, i)$ for $i \in \N$,
where the second component~$i$ is referred to an individual's level.
Each level is occupied by exactly one individual.
The central object in this construction is the \emph{progeny process} $\rho^{\bx} = (\rho^{\bx}_t:\,t \geq 0)$, 
where $\rho_t^{\bx}=(\rho^{\bx}_t(k): k \in \N)$ and for each level~$k$,
$\rho^{\bx}_t(k)$ takes values in $[0,1] \times \Delta_{\dtype}$ (thus carries information about the individual occupying level~$k$ at time~$t$).
For fixed~$t$ and~$k$,
$\rho^{\bx}_t(k)$ codes 
the clonal ancestor of the individual on level~$k$ at time~$t$
(with new clones arising at every mutation event).

The progeny process can be constructed via a lookdown graphical construction. 
We first collect the relevant probabilistic objects.
Let $(\Omega,\mathscr{F},\P)$ be a probability space satisfying the usual conditions and carrying the following independent Poisson point processes:
\begin{itemize}
\item [(I1)] \label{it:I1} $\mathcal{P}_{k,l}$, a Poisson point process on $\R_+$, with intensity measure $\Lambda(\{0\})\dd t$, for every $k<l$ with $k,l\in \N$,
\item [(I2)] \label{it:I2}$\mathcal{N}$, a Poisson point process  on $\R_+\times(0,1]\times [0,1]^{\N}$ with intensity measure $\dd t\times r^{-2}\Lambda(\dd r)\times (\dd u)^{\infty}$,
\item [(I3)] \label{it:I3}$\mathcal{M}_k$, a Poisson point process on $\R_+\times[0,1]$ with intensity measure $\theta\dd t\times \dd u$ for every $k\in\N$.
\end{itemize}

We define the natural filtration $\mathbb{F}:=(\mathscr{F}_t:\, t\geq 0)$,
where for any $t\geq 0$ 
 \begin{equation}\label{eq:Poifiltration}
   \mathscr{F}_t\coloneq\sigma(\mathcal{P}_{k,l}(s),\mathcal{N}(s),\mathcal{M}_k(s):s\in [0,t],\  k<l\in\N).   
 \end{equation}

\begin{figure}[t]
    \centering
    \begin{subfigure}[t]{0.45\linewidth}
        \centering
        \includegraphics[width=1\linewidth]{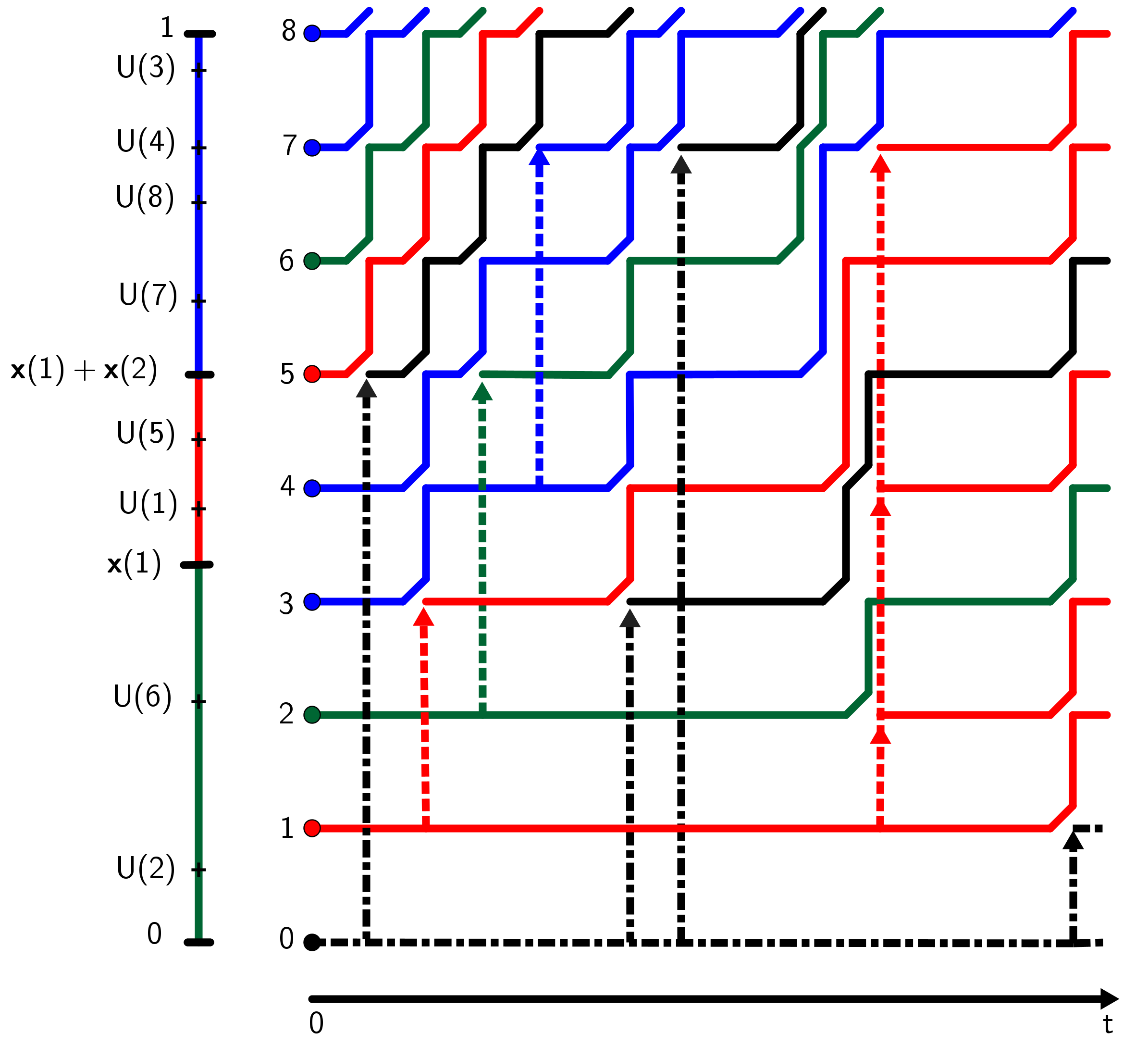}
        \caption{Green, red, and blue lines represent type $1$, type $2$, and type~$3$, respectively. The ghost individual, in dotted, on position~$0$.}
        \label{fig:lookdown_example}
    \end{subfigure}
    \hfill
    \begin{subfigure}[t]{0.45\linewidth}
        \centering
        \includegraphics[width=1\linewidth]{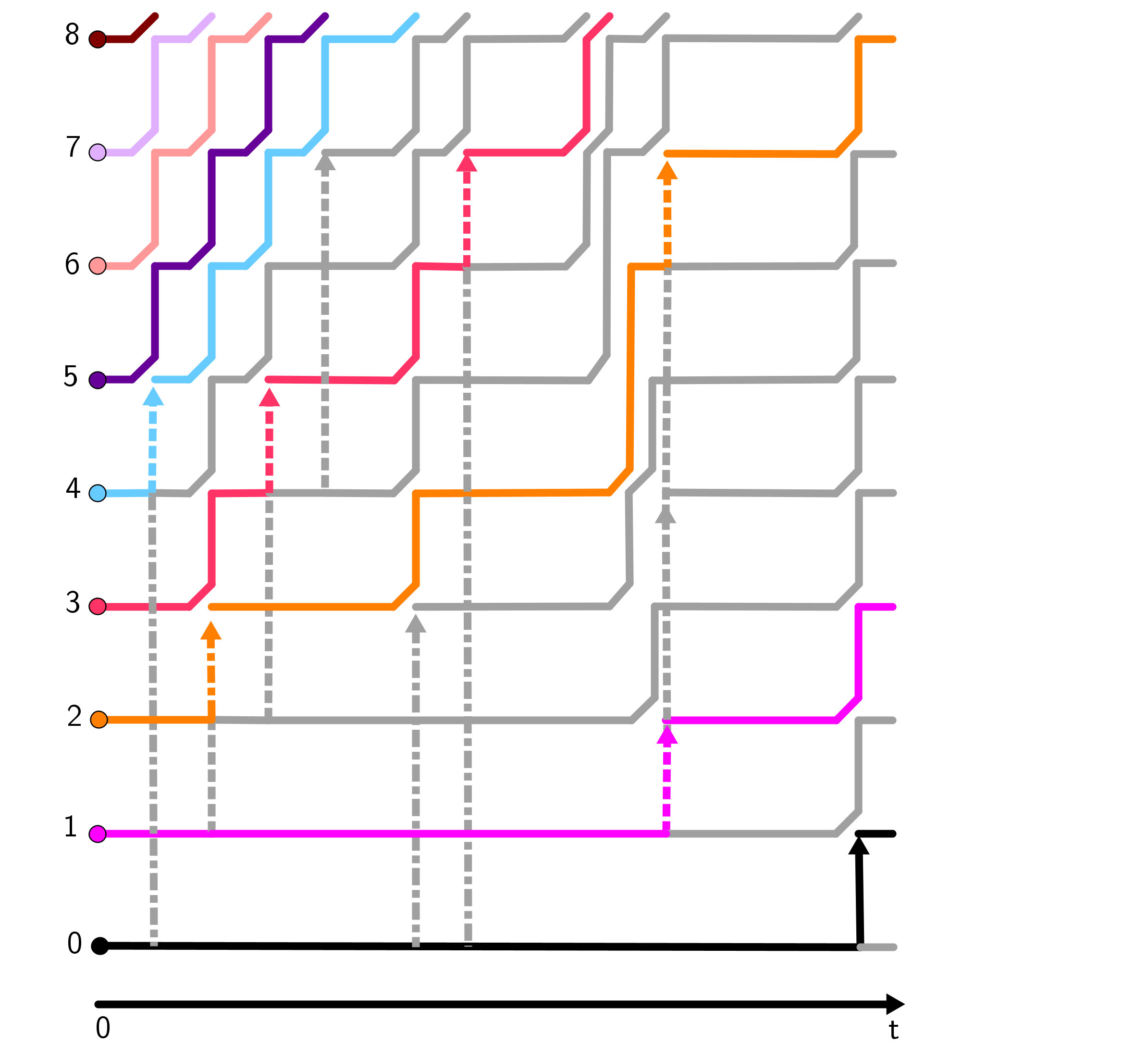}
        \caption{Each color corresponds to a different fixation line in the lookdown-construction pictured on the left.}
        \label{fig:fixation_line}
    \end{subfigure}
    \caption{Example of a lookdown construction and the embedded fixation lines.}
    \label{fig:lookdown}
\end{figure}

Moreover, we require the aforementioned sequence $(\bU(k):\, k\in \N)$ of independent random variables, each uniformly distributed on $[0,1]$, defined on $(\Omega,\mathscr{F},\P)$ and independent of the filtration~$\mathbb{F}$.
Finally, we construct $\rho^{\bx}$ on the basis of the just introduced objects.
Set $\lbl{k}{0}{\bx}:=(\bU(k),\bx)$, 
so the initial individual identifiers are just uniforms,
and their types are determined by the initial distribution~$\bx$.
For $\bx\in \Delta_{\dtype}$ and $k\in \N$, $\rho^{\bx}_\cdot(k):\R_+\to[0,1]\times\Delta_{\dtype}$ is piecewise constant.
We now explain the transitions of the progeny process, see also Fig.~\ref{fig:lookdown_example}.

\begin{figure}[t]
    \centering
    \begin{subfigure}[t]{0.3\linewidth}
        \centering
        \includegraphics[width=.9\linewidth]{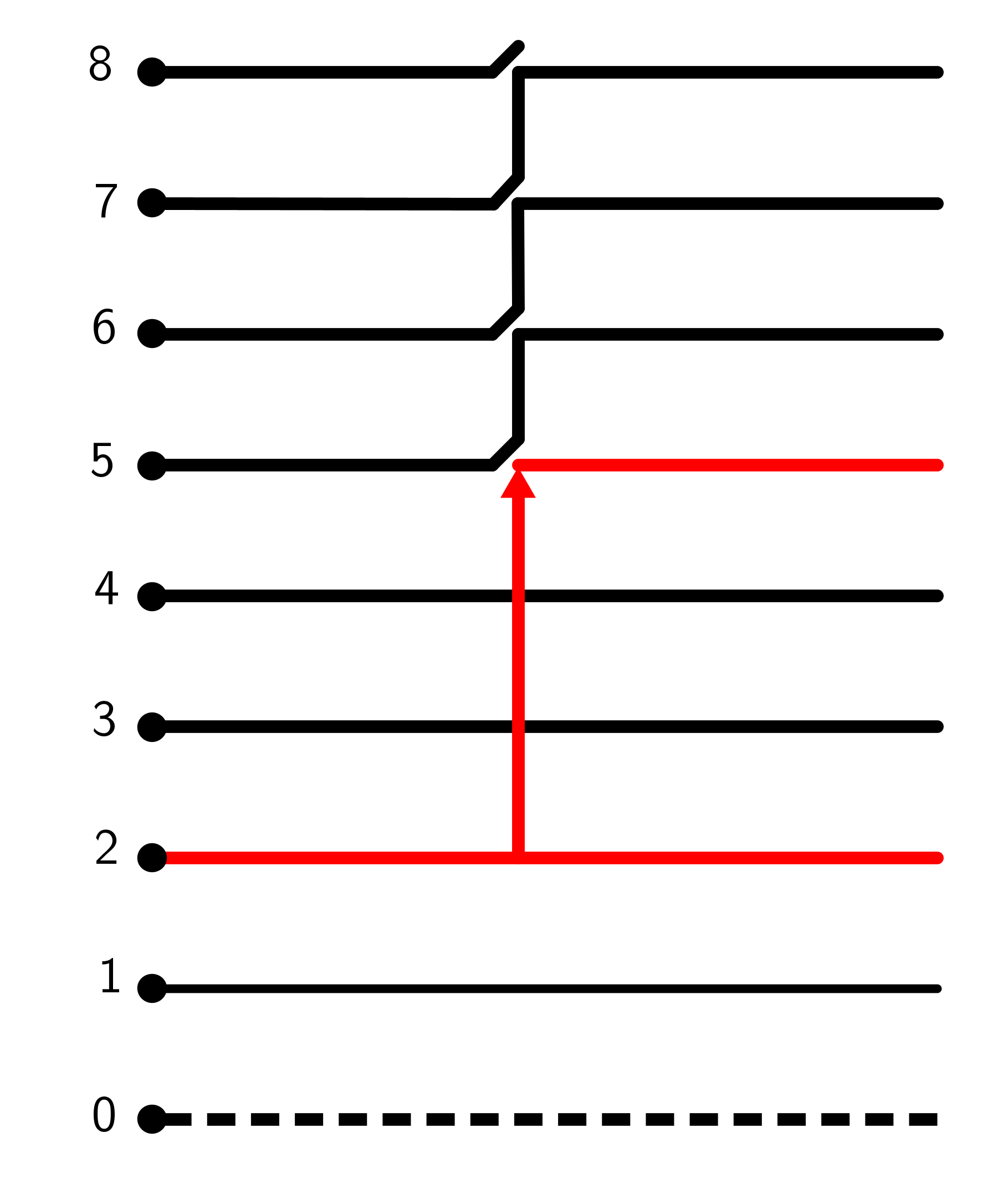}
        \caption{Binary reproduction triggered by a jump of $\mathcal{P}_{2,5}$.}
        \label{fig:binary_reproduction}
    \end{subfigure}
    \hfill
    \begin{subfigure}[t]{0.3\linewidth}
        \centering
        \includegraphics[width=.9\linewidth]{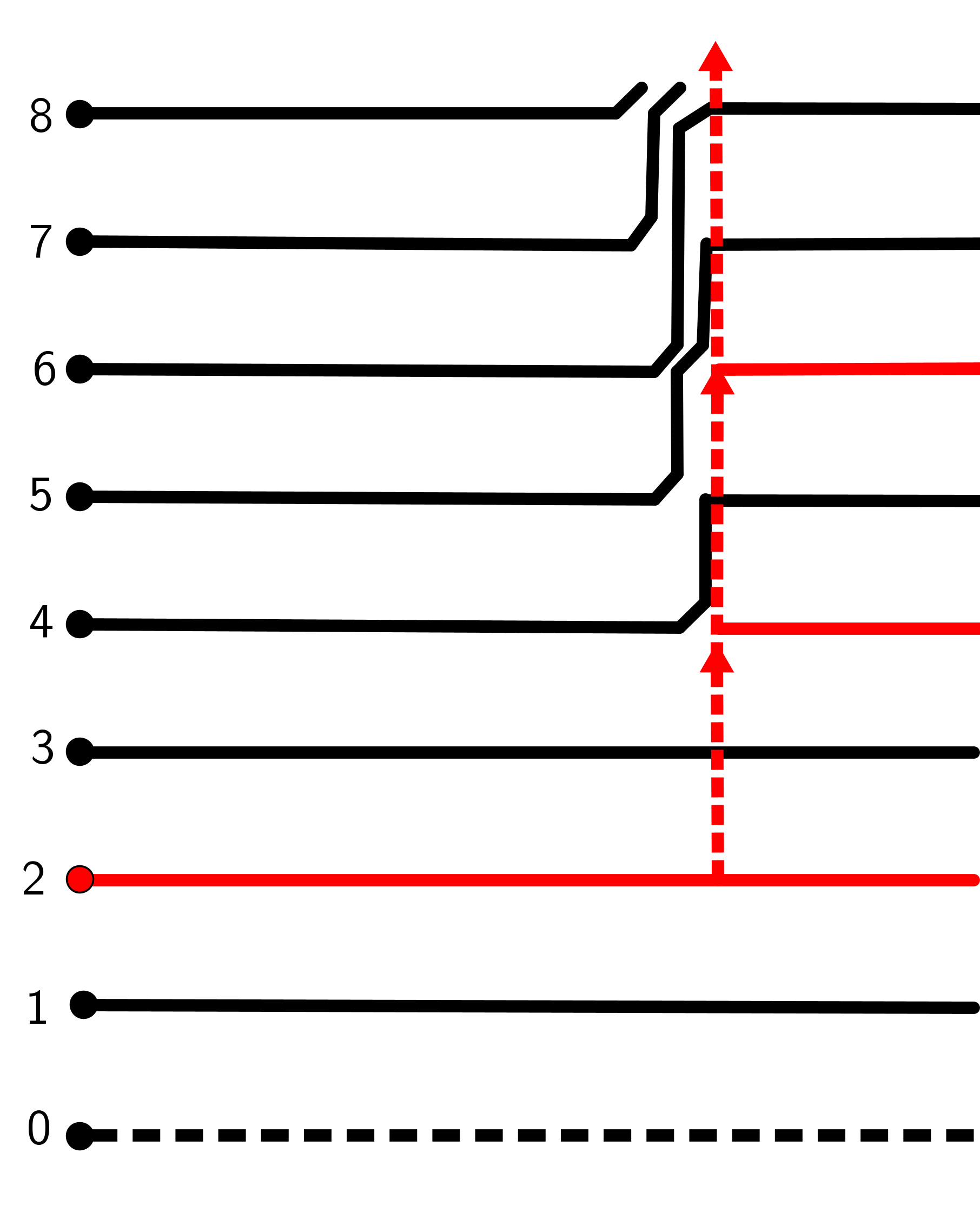}
        \caption{Large reproduction triggered by a jump of $\mathcal{N}$ such that $2$, $4$, $6$ and higher levels are marked.}
        \label{fig:large_reproduction}
    \end{subfigure}
    \hfill
    \begin{subfigure}[t]{0.3\linewidth}
        \centering
        \includegraphics[width=.9\linewidth]{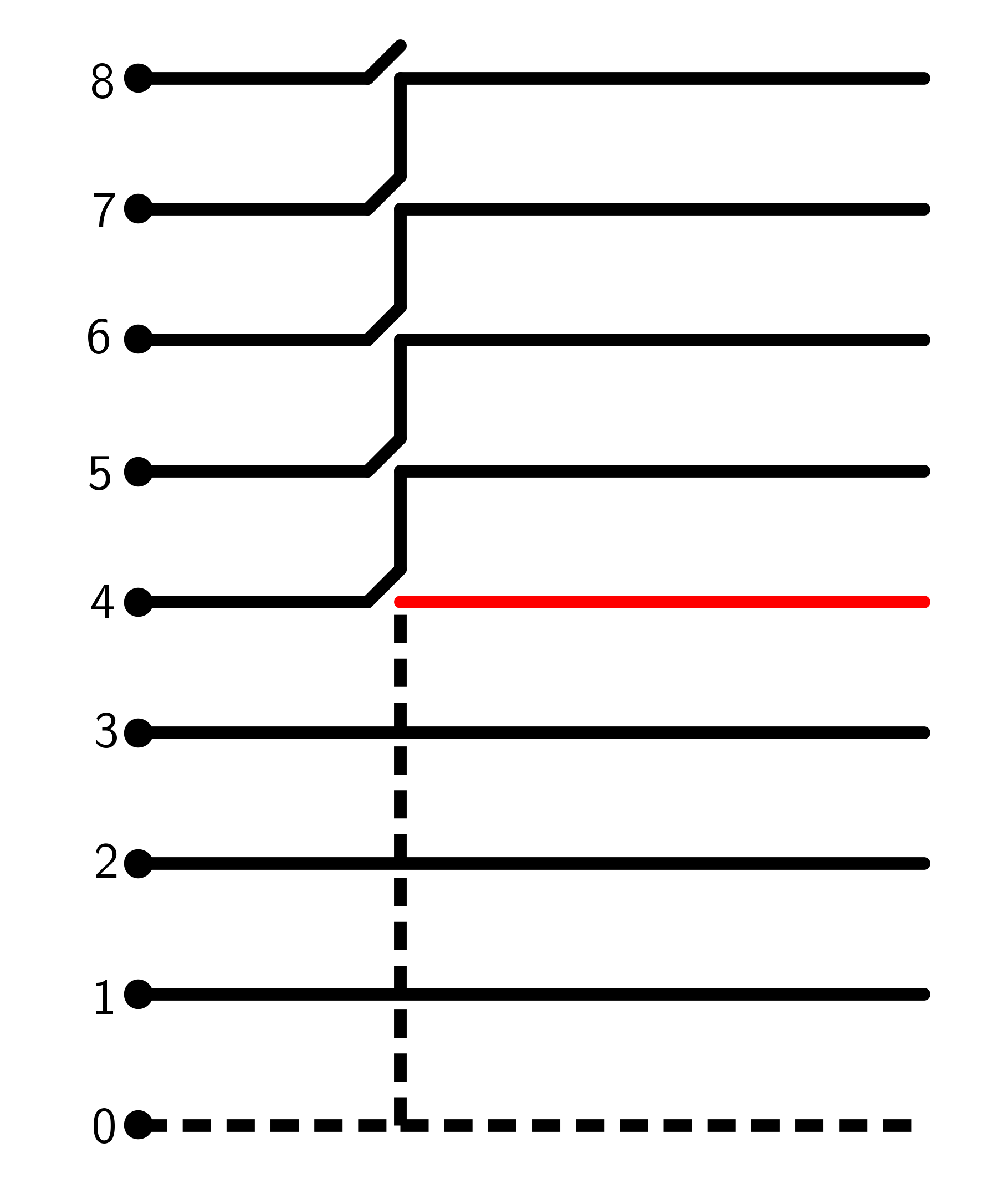}
        \caption{ Mutation triggered by a jump of $\mathcal{M}_4$.}
        \label{fig:mutation}
    \end{subfigure}
    \caption{Events in the lookdown construction.}
    \label{fig:lookdown_transitions}
\end{figure}

\begin{itemize}
    \item [(a)]\textbf{Binary reproduction from $k$ to $l$.} 
    For each atom of $\mathcal{P}_{k,l}$, say at time~$t$, we associate an arrow originating from $(t,k) \in \R_+\times \N$ and pointing to $(t,l)$.
    This arrow corresponds to the event that the individual at level $k$ at time $t$, places a clone of itself at level $l$ at time~$t$.
    Levels that were equal to or above $l$ at time $t-$ are shifted one level up, see also Fig.~\ref{fig:binary_reproduction}.
    In particular, if $t\in \mathcal{P}_{k,l}$, 
    \[\lbl{n}{t}{\bx}=\begin{cases}
		\lbl{n}{t-}{\bx}, &\text{if }n<l,\\
		\lbl{k}{t-}{\bx}, &n=l,\\
		\lbl{n-1}{t-}{\bx}, &n>l.
	\end{cases}
    \]
\end{itemize}

\begin{itemize}
\item  [(b)] \textbf{Large reproduction.} For each atom of $\mathcal{N}$, say $(t, r, \bu)$, proceed as follows. 
At time~$t$, mark each level $l \in \N$ for which $\bu(l) \leq r$. 
Let $k_0$ denote the individual at the smallest marked level, that is, $k_0\coloneqq\min\{l:\, \bu(l)\leq r\}$. 
Then, for each marked level~$l>k_0$, put an arrow from $(t, k_0)$ to $(t, l)$.
This indicates that at time~$t-$, the individual at level~$k_0$ places a clone of itself at all other marked levels.
Individuals are shifted some levels up to make space for the clones, see also Fig.~\ref{fig:large_reproduction}
In particular, if $(t,r,\bu)\in \mathcal{N}$, 
	
	\[\lbl{n}{t}{\bx}=\begin{cases}
		\lbl{n}{t-}{\bx}, &\text{if }n\leq k_0,\\
		\lbl{n-\lvert\{l\leq n:\, \bu(l)\leq r \}\rvert+1 }{t-}{\bx}, &\text{if }\bu(n)>r\text{ and }n>k_0,\\
		\lbl{k_0}{t-}{\bx}, &\text{if }\bu(n)\leq r.
	\end{cases}\]

\item  [(c)] \textbf{Mutation at level $k$.} Finally, for an atom $(t,u)$ of~$\mathcal{M}_k$, 
individuals on levels $k\geq i$ just before time $t-$ move to level $k+1$, 
and a new clone $(u,\bnu)$ is inserted at time~$t$ on position $k$. 
It is instructive to think of the latter as an offspring event on level $k$ from a \emph{ghost} individual on level~$0$ (with each offspring being a new clone), see also Fig.~\ref{fig:mutation}.
In particular, if $(t,u)\in \mathcal{M}_k$, set

 \begin{equation*}\label{eq:mutation}\lbl{n}{t}{\bx}=\begin{cases}
		\lbl{n}{t-}{\bx}, &\text{if }n<k,\\
		(u,\bnu), &n=k,\\
		\lbl{n-1}{t-}{\bx}, &n>k,
	\end{cases}
\end{equation*}
with $\bnu$ being the mutation vector.
\end{itemize}

\begin{lemma}
    The progeny process $\rho^{\bx} = (\rho^{\bx}_t:\,t \geq 0)$ is well-defined.
\end{lemma}
\begin{proof}
The dynamic of $\rho^{\bx}_\cdot(n)$ is only affected by the Poisson point processes $\mathcal{P}_{k,l}, \mathcal{N},\mathcal{M}_k$  with $ k< l\leq n$.
For any arbitrary finite interval $[0,T]$, these processes have a finite number of atoms.
Between their atoms $\rho^{\bx}(k)$ stays constant.
Moreover, the process can be constructed recursively on the levels.
The first level is such that $\rho^{\bx}(1)\coloneqq(\rho^{\bx}_t(1): \, t\geq 0)$ only changes in the atoms of $\mathcal{M}_1$, that is, if $(t,w)\in \mathcal{M}_1$, then set $\rho_t^{\bx}(1)=(w,\bnu)$ and stays in this value until the next atom appears.
If $(\rho_t^{\bx}(k):\, t\geq 0)$ has been constructed for all $k\leq l$, 
then by the above observation $(\rho_t^{\bx}(l+1):\, t\geq 0)$ can be constructed using the levels~$k\leq l$ and the Poisson processes.
    
\end{proof}

\subsection{Convergence}

For each $N \in \N$ and initial type distribution $\bx \in \Delta_{\dtype}$, we define the empirical  types distribution among the first $N$ levels at time $t \geq 0$ by:
\begin{equation*}
	\bX^{\bx,N}_t(i)\coloneqq \frac{1}{N}\sum_{k=1}^N \ind{\mathsf{g}(\lbl{k}{t}{\bx})=i},\qquad i\in[\dtype+1]. 
\end{equation*}
We write $\bX^{\bx,N} = (\bX_t^{\bx,N}:\, t \geq 0)$ for the corresponding stochastic process.

The following result provides a rigorous justification for the lookdown construction as a tool to study the process $\bX^\bx$.
It shows that the empirical distribution at finite levels are in the limit of large~$N$ approximating the distribution of~$\bX^\bx$.
To state the convergence formally, we work in the space  $\mathcal{D}_{\Delta_{\dtype}}([0,\infty))$ of c{\`a}dl{\`a}g-functions with values in $\Delta_{\dtype}$, equipped with the Skorokhod J1-topology;
we denote this convergence by $\xRightarrow{}$.

\begin{theorem}[Convergence] \label{thm:convergence}
	For all $\bx\in \Delta_{\dtype}$,  
	\begin{equation*}
		\bX^{\bx,N} \xRightarrow[N\to\infty]{} \bX^\bx,
	\end{equation*}
	where $\bX^\bx$ is the Markov process with generator~$\mathcal{L}$ defined in \eqref{eq:wfgenerator} and initial condition $\bX_0^\bx=\bx$. 
\end{theorem}

When $\theta = 0$, Theorem~\ref{thm:convergence} becomes a special case of Birkner et al.~\cite[Thm.~1.1]{Birkner2009} for finitely many types. 
However, when $\theta > 0$, that result does not apply, as our mutation mechanism is implemented differently in the lookdown construction, which is the key difference to the construction in Birkner et al.~\cite{Birkner2009}.
In their framework, when an atom $(t, u) \in \mathcal{M}_k$ occurs, the individual at level $k$ and time $t$ is replaced by a new individual $(U, \bnu)$ for a uniform $U$. 
In contrast, our construction shifts each individual at level $l \geq k$ at time $t-$ to level $l+1$, and inserts at level~$k$ a new individual $(U, \bnu)$.
We will prove Theorem~\ref{thm:convergence} at the end of this section.

\subsection{Dual process}
Duality is a useful tool to determine features of a focal process from a dual process, 
with the latter being typically simpler to analyze.
For us, the dual process is pivotal in proving Theorem~\ref{thm:convergence}.
In our case, we call the dual process the \emph{ancestral process},
and it is defined on the basis of yet another process that we refer to as \emph{color configuration process}.
This latter process arises when trying to determine in the lookdown construction the probability of a given configuration of types on the levels at a fixed time,
conditionally on the ancestry of the population.

The color configuration process $\bC=(\bC_r:\, r\geq0)$ takes values in the set $(\cup_{k\geq0} [d]^k)\cup \dagger$, 
where~$\dagger$ is a cemetery state. 
For $\bc\in (\cup_{k\geq0} [d]^k)$, we write $\lvert \bc\rvert \coloneqq \dim(\bc)$.
The construction relies on the Poisson point processes introduced earlier. 
To explain this, it is instructive to fix $t>0$ and consider $(\bC_r:\, r\in[0,t])$,
because then the relation to the progeny process is most direct.
Let $\bC_0=\bc_0\in\cup_{k\geq0} [d]^k$.
We view an initial color process configuration as $\lvert \bC_0 \rvert$ colored balls (or ancestral types) in the lookdown construction at time~$t$; so $\bC_0(l)$ should be interpret as the color of the ball placed on level~$l$ at time~$t$ in the lookdown construction.
The process $\bC = (\bC_r:\, r \in [0,t])$ then evolves backward in time (from time~$t$ to time~$0$) by following the arrows in the graphical representation.

\begin{figure}[t]
    \centering
    \begin{subfigure}[t]{0.47\linewidth}
        \centering
        \includegraphics[width=.9\linewidth]{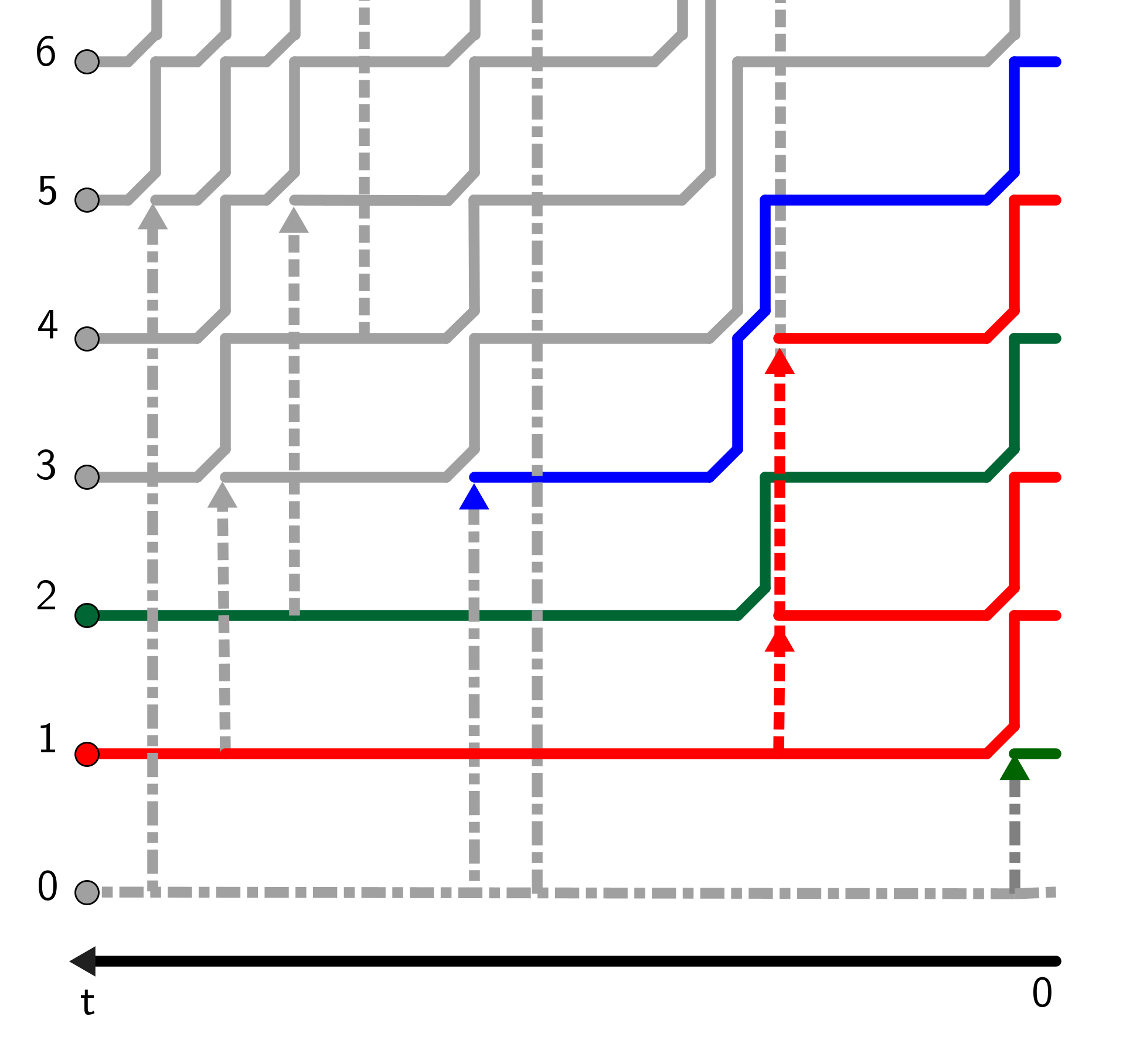}
        \caption{ Realization of the color configuration process without unfeasible events. The first event is a mutation at level~$1$, followed by a large reproduction involving levels $1$, $2$, and $4$, and a mutation on level~$3$.}
        \label{fig:without_unfeasible}
    \end{subfigure}
    \hfill
    \begin{subfigure}[t]{0.47\linewidth}
        \centering
        \includegraphics[width=.9\linewidth]{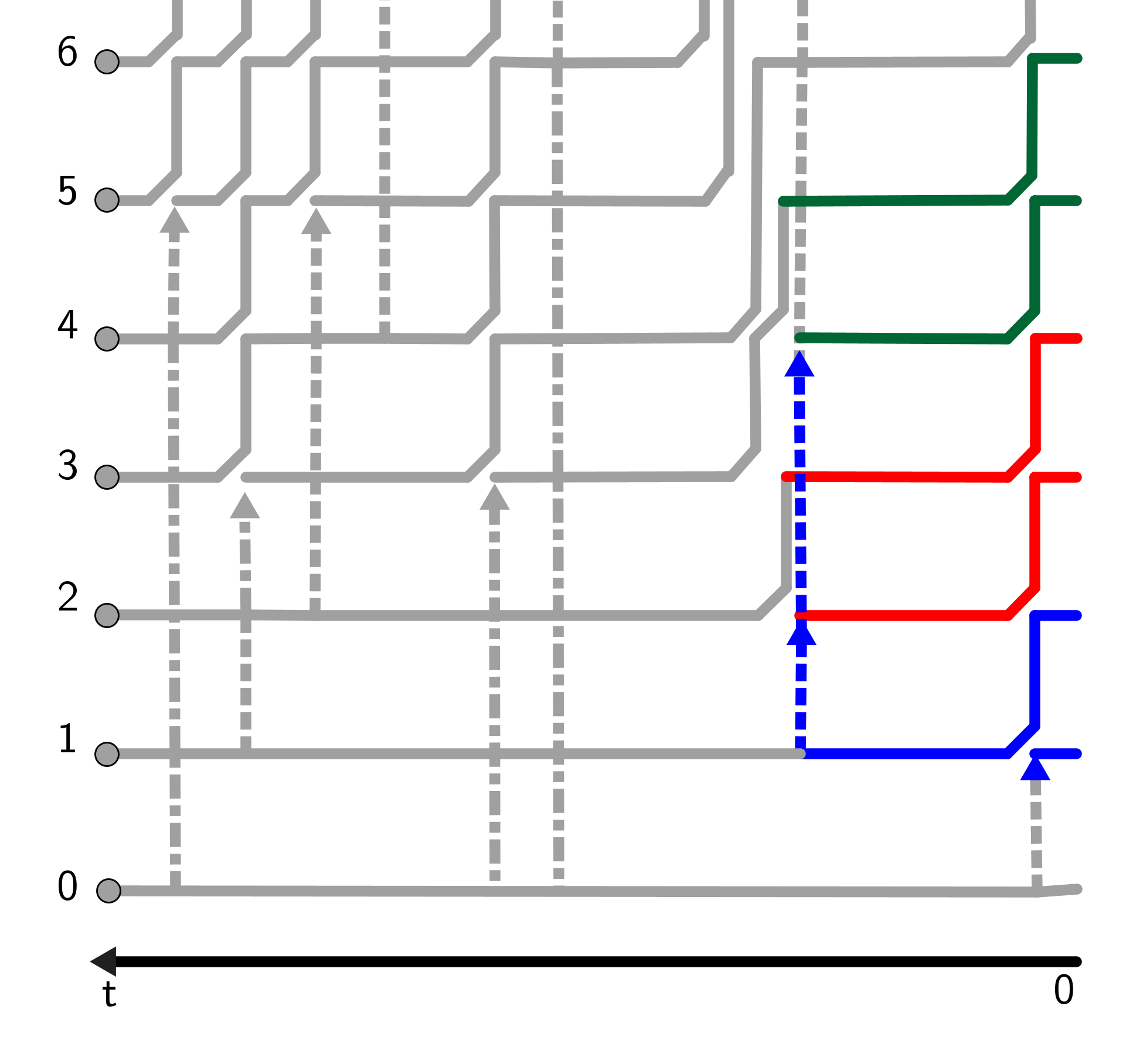}
        \caption{Realization of the color configuration process with an unfeasible event.
        The first event is a mutation on level~$1$.
        The second event is a large reproduction on level $1$, $2$ and $4$.
        Because the levels have different colors, the event is unfeasible and the process transitions to~$\dagger$.}
        \label{fig:with_unfeasible}
    \end{subfigure}
    \caption{Examples of feasible and unfeasible events in the color configuration process.}
    \label{fig:colorprocess}
\end{figure}

If the color configuration process is in~$\dagger$, the process absorbs there and does not change anymore.
Otherwise, the configuration process may change only at the jump times of the Poisson point processes.
At such times, we consider the balls with levels associated to the arrows in the lookdown construction:
if the balls have different colors, the process is sent to the absorbing (cemetery) state~$\dagger$, see cases~(i) below. 
Otherwise, if all balls involved in an event have the same color, some of them may be removed, see cases~(ii) below. 
The precise update rule at an atom of the Poisson point processes are as follows. 
See also Fig.~\ref{fig:colorprocess} for an illustration of the process.

\begin{itemize}
	\item  [(a)] \textbf{Binary reproduction.} For $t-r\in \mathcal{P}_{k,l}$ with $1\leq k<l\leq \lvert \bC_{r-}\rvert$:
	\begin{enumerate}
		\item[(i)] \textbf{Unfeasible event.} If $\bC_{r-}(k)\neq \bC_{r-}(l)$, set $\bC_r=\dagger$.  
        
		\item[(ii)]\textbf{Kingman's coalescence.} Otherwise set $\lvert\bC_r\rvert=\lvert\bC_{r-}\rvert-1$ and 
	\[\bC_r(m)=\begin{cases}
		\bC_{r-}(m), &\text{if }1\leq m<l,\\
		\bC_{r-}(m+1), &\text{if }l\leq m \leq \lvert\bC_{r-}\rvert-1.
	\end{cases}\]
	\end{enumerate}
    
	\item  [(b)] \textbf{Large reproduction.} For $(t-r,p,\bu)\in \mathcal{N}$ and $\lvert\{l\in[\lvert\bC_{r-}\rvert ]: \bu(l)\leq p\}\rvert>1$:
	\begin{enumerate}
		\item[(i)] \textbf{Unfeasible event.} If $\lvert\{\bC_{r-}(l):\bu(l)\leq p\}\rvert>1$, set $\bC_r=\dagger$.
        
		\item[(ii)] \textbf{$\Lambda$-coalescence.} Otherwise set $\lvert\bC_r\rvert=\lvert\bC_{r-}\rvert-\lvert\{l\in[\lvert\bC_{r-}\rvert]: \bu(l)\leq p\}\rvert+1$ and remove the balls on all level $l$ with $\bu(l)\leq p$
        except the first of such balls.
        The levels are updated as follows: for any $m\leq \lvert\bC_r\rvert$,
		\[\bC_r(m)=\begin{cases}
			\bC_{r-}(m), &\text{if } \lvert \{l\leq m: \bu(l)\leq p\}\rvert\leq 1,\\
			\bC_{r-}(m+\lvert \{l\leq m: \bu(l)\leq p\}\rvert-1), &\text{otherwise}.
		\end{cases}\]
	\end{enumerate}
    
	\item  [(c)] \textbf{Mutation at level $\mathbf{k}$.} For $(t-r,u)\in \mathcal{M}_k$ with $k\leq \lvert\bC_{r-}\rvert$:
	\begin{enumerate}
		\item[(i)] \textbf{Unfeasible event.}  If $\bC_{r-}(k)\neq \type{u}{\bnu}$, set $\bC_r=\dagger$.
        
		\item[(ii)] \textbf{Mutation.} Otherwise set  $\lvert\bC_r\rvert=\lvert\bC_{r-}\rvert-1$, remove ball $k$ and update the levels, that is
		\[\bC_r(l)=\begin{cases}
			\bC_{r-}(l), &\text{if } l<k,\\
			\bC_{r-}(l+1), &\text{if }k\leq l \leq \lvert\bC_{r-}\rvert-1.
		\end{cases}\]
	\end{enumerate}
\end{itemize}

The \emph{ancestral process} is then defined on the basis of the color configuration process started from a particular initial configuration.
The state space is $\N_{0,\dagger}^{\dtype}\coloneqq \N_0^{\dtype}\cup{\dagger}$.
For an initial state $\bn\in \N_0^d$ of the ancestral process, 
we consider a color configuration process started in $\bc^{(0)}(\bn)\in \cup_{l\geq0} [d]^l$, 
where the first $\bn(1)$ entries are~$1$, 
the next $\bn(2)$ entries are $2$, 
$\ldots$, and the last $\bn(\dtype)$ entries are $\dtype$.
So for our interpretation this means that at time~$t$ in the lookdown construction
the first $\bn(1)$ balls are type (or color) $1$ and they are placed on the first $\bn(1)$ levels, 
the next  $\bn(2)$ balls are type (or color) $2$ and they are placed on the next $\bn(2)$ levels, and so on, 
with the last $\bn( d)$ balls being of type (or color) $d$.
For every fixed $t>0$ the ancestral process $\bA^{(\bn)}=(\bA_r^{(\bn)}:\, r\in [0,t])$ starts in $\bn$.
Let $(\bC_r:\, r\in [0,t])$ be the corresponding color configuration process in~$[0,t]$ (started at time~$t$ from $\bc^{(0)}(\bn)$).
For $r\in [0,t]$, if $\bC_r = \dagger$, then $\bA_r^{(\bn)}=\bC_r=\dagger$;
otherwise, if $\bC_r \neq \dagger$,
then we have $\bA^{(\bn)}_r = (\bA^{(\bn)}_r(i):\, i\in [d])$ with 
\[
\bA^{(\bn)}_r(i):=\lvert\{l:\bC_{r}(l)=i\}\rvert, \qquad i\in [d],
\] 
i.e. $\bA^{(\bn)}_r(i)$ gives the number of balls in $\bC_r$ with color $i$.

Thus, $\bA^{(\bn)}=(\bA_r^{(\bn)}:\, r\in [0,t] )$ is a  Markov chain taking values in $\N_{0,\dagger}^{\dtype}$, starting at $\bn$ 
and with infinitesimal generator
\begin{equation}\label{eq:ancestralgen}
\mathcal{A}g(\bn)=\mathcal{A}_{\mathrm{WF}}g(\bn)+\mathcal{A}_{\Lambda}g(\bn)+\mathcal{A}_{\theta}g(\bn), 
\qquad g\in \mathcal{B}(\N_{0,\dagger}^{\dtype})   
\end{equation}
where
\begin{align*}
    \mathcal{A}_{\mathrm{WF}}g(\bn)
    &=\Lambda(\{0\})\sum_{i=1}^d \bigg(\bn(i)(\bn(i)-1)[g(\bn-\e{i})-g(\bn)]+\sum_{\substack{j=1 \\ j\neq i}}^d \bn(i)\bn(j)[g(\dagger)-g(\bn)]\bigg), \\
    \mathcal{A}_{\Lambda}g(\bn)
    &=\sum_{i=1}^d\sum_{k=2}^{\bn(i)}\binom{\bn(i)}{k}\lambda_{\lvert \bn \rvert,k} [g(\bn-(k-1)\e{i})-g(\bn)]		\\
    &\quad + \int_{(0,1]} \bigg(1-(1-r)^{\lvert \bn\rvert} -\sum_{i=1}^d\sum_{k=1}^{\bn(i)} \binom{\bn(i)}{k} (1-r)^{\lvert \bn\rvert-k} r^{k}\bigg)\frac{\Lambda(\dd r)}{r^2} [g(\dagger)-g(\bn)], \\
    \mathcal{A}_{\theta}g(\bn)
    &=\theta\sum_{i=1}^d \bn(i)\bigg( \bnu(i)[g(\bn-\e{i})-g(\bn)]+\sum_{\substack{j=1\\ j\neq i}}^d \bnu(j)[g(\dagger)-g(\bn)]\bigg),
\end{align*}
$\mathcal{A}_{\mathrm{WF}}g(\dagger)=	\mathcal{A}_{\Lambda}g(\dagger)=	\mathcal{A}_{\theta}g(\dagger)=0$, and $\lambda_{n,k}$ the transition rate of a $\Lambda$-coalescent defined in \eqref{eq:ratetransLam}.

We remark that for $\mathcal{A}_{\Lambda}$, the rate at which a transition to the absorbing state~$\dagger$ occurs can be rewritten as 
\[
\int_{(0,1]} \bigg( 1-(1-r)^{\lvert \bn\rvert} - \lvert \bn\rvert (1-r)^{\lvert \bn\rvert-1}r-
\sum_{i=1}^d\sum_{k=2}^{\bn(i)} \binom{\bn(i)}{k} (1-r)^{\lvert \bn\rvert-k} r^{k}\bigg)\frac{\Lambda(\dd r)}{r^2},
\]
and thus corresponds to the rate at which at least two levels are selected to coalesce, and these levels are not of the same type.

Our next result establishes the duality between $\bX^{\bx}$ and $\bA^{(\bn)}$ via the function  
$H:\Delta_{d}\times\N_{0,\dagger}^{\dtype}\to [0,1]$ defined for $\bx \in\Delta_{d}$ and $\bn\neq\dagger$ by
\begin{equation}\label{definition H}
	H(\bx,\bn)=\bx^{\bn}:= \prod_{k=1}^{d}\bx(k)^{\bn(k)},
\end{equation}
and $H(\bx,\dagger)=\bx^{\dagger}\coloneqq0$. 

\begin{theorem}[Moment duality] \label{thm:duality}
For $\bx\in \Delta_{\dtype}$ and $\bn\in \N_{0,\dagger}^{\dtype}$, 
let $\bX^{\bx}$ and $\bA^{(\bn)}$ be the processes governed by generator $\mathcal{L}$ and $\mathcal{A}$, respectively, and with respective initial state $\bx$ and $\bn$.
Then, for all $t\geq 0$,
 \[\E[H(\bX_t^{\bx},\bn)]=\E[H(\bx,\bA_t^{(\bn)})].\]
\end{theorem}
\begin{proof}
    For any $\bn$, the process $\bA^{(\bn)}$ is a Markov chain with finite state space. 
    Additionally, $\bX^{\bx}$ is a Feller process because the one-dimensional projections $(\bX^{\bx}_t(i):\, t\geq 0)$, $i\in [d]$,
    are a $2$-type Fleming-Viot processes which are known to have the Feller property, see for instance~\cite[Prop.~B.5]{Cordero2022}. 
    Then, by Liggett~\cite[Thm. 3.42]{liggett2010}, the result holds once  we prove
    \begin{equation}
        \mathcal{L}H(\bx,\bn)=\mathcal{A}H(\bx,\bn),\qquad \mbox{for all } \bx,\bn.\label{eq:dualityidentity}
    \end{equation} 
    Hence, we are going to compare the respective parts of the generators \eqref{eq:wfgenerator} and \eqref{eq:ancestralgen}. 
    First we note that
\[
        \frac{\partial}{\partial\bx(i)} H(\bx,\bn)=\bn(i)H(\bx,\bn-\e{i})\]
and
\[ \frac{\partial^2}{\partial\bx(i)\partial\bx(j)} H(\bx,\bn)=\bn(i)(\bn(j)-\ind{i}(j))H(\bx,\bn-\e{i}-\e{j}).
\]
    Using this, we compute
    \begin{align*}
        \mathcal{L}_{\mathrm{WF}}H(\cdot,\bn)(\bx)&=\Lambda(\{0\})\sum_{i=1}^d\Bigg(\binom{\bn(i)}{2} [H(\bx,\bn-\e{i})-H(\bx,\bn)]-\sum_{\substack{j=1\\ j\neq i}}^d \bn(i)\bn(j)H(\bx,\bn) \Bigg)\\
        &=\mathcal{A}_{\mathrm{WF}}H(\bx,\cdot)(\bn),
    \end{align*}
    and 
    \begin{align*}
        \mathcal{L}_{\theta}H(\cdot,\bn)(\bx)&=\theta\sum_{i=1}^d\bn(i)\Bigg(\bnu(i)[H(\bx,\bn-\e{i})-H(\bx,\bn)]+\sum_{j\neq i}\bnu(j)[H(\bx,\dagger)-H(\bx,\bn)] \Bigg)\\
        &=\mathcal{A}_{\theta}H(\bx,\cdot)(\bn).
    \end{align*}
    
    For the $\Lambda$-part, the computation is slightly more involved:
    \begin{align*}
        \mathcal{L}_{\Lambda}H(\cdot,\bn)(\bx) 
        =&\int_{(0,1]} \Bigg(\sum_{i=1}^{\dtype +1} H((1-r)\bx+r\e{i},\bn)\Bigg)-H(\bx,\bn)\frac{\Lambda(\dd r)}{r^2}\\
        =&\int_{(0,1]} \Bigg( \Bigg\{ \sum_{i=1}^d \Bigg[\sum_{k=0}^{\bn(i)} \binom{\bn(i)}{k} r^{k}\bx(i)^{\bn(i)-k+1}(1-r)^{\lvert \bn \rvert-k} \prod_{\substack{\ell=1\\ \ell \neq i}}^d \bx(\ell)^{\bn(\ell)}\Bigg]\\
        &+ \bigg(1-\sum_{i=1}^d \bx(i) \bigg) (1-r)^{\vert \bn\rvert} \bx^{\bn}\Bigg\} - \bx^{\bn} \Bigg)\frac{\Lambda(\dd r)}{r^2}
    \end{align*}
    Then, by arranging the terms 
    	\begin{align*}
    		 \mathcal{L}_{\Lambda}H(\cdot,\bn)&(\bx) 
    		=\int_{(0,1]} \Bigg(  \sum_{i=1}^d \Bigg[\sum_{k=2}^{\bn(i)} \binom{\bn(i)}{k} r^{k}\bx(i)^{\bn(i)-k+1} (1-r)^{\lvert \bn \rvert-k} \prod_{\substack{\ell=1\\ \ell \neq i}}^d \bx(\ell)^{\bn(\ell)}\Bigg] \\
    		&- \Big(1-(1-r)^{\lvert \bn\rvert} - \lvert \bn\rvert r(1-r)^{\lvert \bn\rvert -1}\Big) \bx^{\bn} \Bigg)\frac{\Lambda(\dd r)}{r^2} \\
    		=&\sum_{i=1}^d \sum_{k=2}^{\bn(i)} \binom{\bn(i)}{k}\lambda_{\lvert \bn\rvert ,k}\Big[H(\bx,\bn-(k-1)\e{i})-H(\bx,\bn)\Big] + [H(\bx,\dagger)-H(\bx,\bn)]\\
    		&\cdot  \int_{(0,1]}\Bigg(1-(1-r)^{\lvert \bn\rvert} - \lvert \bn\rvert r(1-r)^{\lvert \bn\rvert -1} -\sum_{i=1}^d\sum_{k=2}^{\bn(i)} \binom{\bn(i)}{k} r^k(1-r)^{\lvert \bn\rvert -k} \Bigg) \frac{\Lambda(\dd r)}{r^2}\\
    		=&\mathcal{A}_{\Lambda}H(\bx,\cdot)(\bn).
    	\end{align*}
    Therefore, each building block of one generator satisfies the equality with the respective part of the other one,
    and thus \eqref{eq:dualityidentity} holds.
\end{proof}

\subsection{Proof of Theorem~\ref{thm:convergence}} 
The duality links the process $\bX^{\bx}$ with the ancestral process $\bA^{(\bn)}$. In order to prove Theorem~\ref{thm:convergence} we will need two lemmas. In them, we connect the ancestral process both with the progeny process $\lbl{\cdot}{t}{\bx}$ and with the limit of the empirical type distributions $\bX^{\bx,N}$. For this purpose we define, for any $\bn \in \N_0^d$, the event $E_{\bx,\bn,t}$. This is the event that, in the lookdown construction at a fixed time $t$, the first $\bn(1)$ levels are of type~$1$, the next $\bn(2)$ levels are of type~$2$, and so on, with the last $\bn(d)$ levels of type~$d$; that is,
\begin{equation}\label{eq:Exnt}
E_{\bx,\bn,t}= \bigcap_{\substack{k\in[d] \\ \bn(k)>0 }}\left\{ \mathsf{g}\left(\lbl{1+\sum_{i=1}^{k-1}\bn(i)}{t}{\bx}\right)=\cdots=\mathsf{g} \left(\lbl{ \bn(k)+ \sum_{i=1}^{k-1}\bn(i)}{t}{\bx}\right)=k\right\}.
\end{equation}
The probability of $E_{\bx,\bn,t}$ can be computed by tracing the ancestral lineages backward in time and requiring that all ancestors of type-$i$ levels are of type~$i$.

\begin{lemma}\label{lem:momentstypes}
	For all $\bx\in \Delta_{d}$, $\bn\in \N_0^d$ and $t\geq 0$, 
	\begin{equation*}
	\P\left( E_{\bx,\bn,t} \right) 
	= \E\left[\bx^{\bA_t^{(\bn)}}\right] .
	\end{equation*}
\end{lemma}

\begin{proof}
We can decompose $E_{\bx,\bn,t}$ into all the possible starting configurations of $(\mathsf{g}(\lbl{k}{0}{\bx}):\, k\in \N)$. More precisely, 
\[
E_{\bx,\bn,t}=\underset{\bg_0\in[d]^{\infty}}{\bigcup} \Big(E_{\bx,\bn,t}\cap \left\{\mathsf{g}(\lbl{k}{0}{\bx})=\bg_0(k):\,  k\in \N\right\}\Big).
\]
First, we  condition on the information of the Poisson processes given by the $\sigma$-algebra, $\mathscr{F}_t$ defined in~\eqref{eq:Poifiltration}.
Observe that the progeny process  $\rho^{\bx}$ and the color configuration process $\bC$ use the same arrows associated with $\mathscr{F}_t$; the first process starts with a random configuration (independent of $\mathscr{F}_t$) and propagates the information by using the arrows forward in time.
The second process starts with a deterministic configuration and propagates the information by using the arrows backward in time, which implies that~$\bC$ is $\mathscr{F}_t$ measurable.
If we condition on $\mathscr{F}_t$, the set where there is not an initial configuration $\bg_0$ such that forward in time the configuration at time $t$ restricted to the first $\lvert \bn \rvert$ levels is equal to $\bC_0$, is equal to the sets $\{E_{\bx,\bn,t}=\emptyset\}$ and $\{\bC_t= \dagger\}$.
In the complement set, 
$E_{\bx,\bn,t}\cap \left\{\mathsf{g}(\lbl{k}{0}{\bx})=\bg_0(k):\, k\in \N\right\}\neq \emptyset$ if and only if 
$\bg_0(k)=\bC_t(k)$ for every $1\leq k \leq  \lvert\bC_{t}\rvert$.
Then,

\begin{equation*}
    \begin{split}
        \P&(E_{\bx,\bn,t}\mid \mathscr{F}_t)\\
       & =\E\left(\ind{\bC_t\neq \dagger}
        \underset{\bg_0\in[d]^{\infty}}{\sum}  \ind{\mathsf{g}(\lbl{k}{0}{\bx})=\bg_0(k)=\bC_t(k):\; k\leq \lvert \bC_{t}\rvert}\ind{\mathsf{g}(\lbl{k}{0}{\bx})=\bg_0(k):\; k> \lvert\bC_{t}\rvert} \,\Big\vert\, \mathscr{F}_t \right)
    \end{split}
\end{equation*}	
Now, we use that $\bC_t$ is $\mathscr{F}_t$ measurable, the variables $(\bU(k): k\in \N)$ are i.i.d.  uniform random variables in $[0,1]$,
independent of $\mathscr{F}_t$, and the definition of $\mathsf{g}(\lbl{k}{0}{\bx})=\type{\bU(k)}{\bx}$, to obtain
\begin{equation*}
    \begin{split}
        \P&(E_{\bx,\bn,t}\mid \mathscr{F}_t)\\
&=\ind{\bC_t\neq \dagger} \underset{\bg_0\in[d]^{\infty}}{\sum}\prod_{k=1}^{\lvert\bC_{t}\rvert}
        \P\left(
         \type{\bU(k)}{\bx}=\bg_0(k)=\bC_t(k)\right) \prod_{k=1+\lvert\bC_{t}\rvert}^{\infty}
        \P\left(\type{\bU(k)}{\bx}=\bg_0(k) \right)\\
        &=\ind{\bC_t\neq \dagger}
        \prod_{k=1}^{\lvert\bC_{t}\rvert}
        \bx(\bC_t(k)).
    \end{split}
\end{equation*}	
By observing that the last equation is equal to $\bx^{\bA_t^{(\bn)}}$ and taking the expectation, we obtain the result.

\end{proof}

Now, we  analyze the limit of the empirical type distributions $\bX^{\bx,N}$. The following lemma establish that the moments of the processes  $\bX^{\bx,N}$ converges to the generating function of $\bA^{(\bn)}$, see Lemma \ref{lem:momentstypes}.

\begin{lemma} \label{lem:asymmomentstypes}
	For all $\bx\in \Delta_d$, $\bn\in \N_{0,\dagger}^{\dtype}$ and $t\geq0$,
	\begin{equation}\label{eq:prob}
		\E[(\bX_t^{\bx,N})^{\bn}]\xrightarrow[N\to\infty]{}\P\left(E_{\bx,\bn,t}\right).
	\end{equation}
\end{lemma}

\begin{proof}
		Fix $\bx\in \Delta_d$, $\bn\in \N_0^{d}$, $N\geq \lvert \bn\rvert$ and $t\geq 0$. 
	Using the definition of $\bX_t^{\bx,N}$ and the multinomial theorem, we obtain	
	\begin{align*}
		\E&\left[\prod_{k=1}^d(\bX_t^{\bx,N}(k))^{\bn(k)}\right]
=\frac{1}{N^{\lvert \bn \rvert}} \E\left[\prod_{k=1}^d\left(\sum_{\substack{i_1^k+\cdots+i_N^k=\bn(k) \\ i_1^k,\ldots,i_N^k\geq 0}}\binom{\bn(k)}{i_1^k,\ldots,i_N^k}\prod_{j=1}^N \ind{\mathsf{g}(\lbl{j}{t}{\bx})=k}^{i_j^k}\right)\right]\\
&\hspace{1.1cm}=\frac{1}{N^{\lvert \bn \rvert}} \sum_{\substack{i_1^1+\cdots+i_N^1=\bn(1) \\ i_1^1,\ldots,i_N^1\geq 0}}\dots\sum_{\substack{i_1^d+\cdots+i_N^d=\bn(d) \\ i_1^d,\ldots,i_N^d\geq 0}}\E\Bigg[\prod_{k=1}^d\binom{\bn(k)}{i_1^k,\ldots,i_N^k}\prod_{j=1}^N \ind{\mathsf{g}(\lbl{j}{t}{\bx})=k}^{i_j^k}\Bigg].
	\end{align*}
Now, we consider  $I=| \{i_j^k: i_j^k\neq 0, j\leq N, k\leq d\}|$ and we decompose the previous sum according to the value of $I$,  
	\begin{align}
&	\E[(\bX_t^{\bx,N})^{\bn}]\label{eq:upperlower} \\
	&=\frac{1}{N^{\lvert \bn \rvert}} \sum_{\ell=0}^{\lvert \bn\rvert}
 \sum_{\substack{i_1^1+\cdots+i_N^1=\bn(1) \\ i_1^1,\ldots,i_N^1\geq 0}}\dots\sum_{\substack{i_1^d+\cdots+i_N^d=\bn(d) \\ i_1^d,\ldots,i_N^d\geq 0}}\ind{ I=\ell}
\E\Bigg[\prod_{k=1}^d\binom{\bn(k)!}{i_1^k,\ldots,i_N^k}\prod_{j=1}^N \ind{\mathsf{g}(\lbl{j}{t}{\bx})=k}^{i_j^k}\Bigg].\nonumber
	\end{align}

We will establish the result by  showing that the term associated with $I=\lvert \bn\rvert$ converges to the probability on the right-hand side of~\eqref{eq:prob} and 
the terms associated with $I<\lvert \bn\rvert$ vanish as $N\to\infty$. 
First, we consider when ~$I=\lvert \bn\rvert$ and note that in this case  $i_j^k\in \{0,1\}$ for any $j\leq N$ and $k\leq d$. Therefore, the term associated with  $I=|\bn|$  in \eqref{eq:upperlower} reduces to 
	\begin{equation}\label{eq:auxlem31}
	\frac{\bn(1)!\cdots \bn(d)!}{N^{\lvert \bn \rvert}} 	\sum_{\substack{i_1^1+\cdots+i_N^1=\bn(1) \\ i_1^1,\ldots,i_N^1\geq 0}}\dots\sum_{\substack{i_1^d+\cdots+i_N^d=\bn(d) \\ i_1^d,\ldots,i_N^d\geq 0}}\ind{ I=\lvert \bn\rvert}
\E\Bigg[\prod_{k=1}^d\prod_{j=1}^N \ind{\mathsf{g}(\lbl{j}{t}{\bx})=k}^{i_j^k}\Bigg].
	\end{equation}
Observe that if $i_j^k>0$ and $\mathsf{g}(\lbl{j}{t}{\bx})=k$ for some $j\leq N$ and $k\leq d$, then, for a positive contribution of the summand it is necessary that $i_j^l=0$ for the other types $l\neq k$ because $\ind{\mathsf{g}(\lbl{j}{t}{\bx})=l}=0$. 	
Since $\lvert\{i_j^k: i_j^k=1\}\rvert=\lvert \bn\rvert$, for a positive contribution we need to select $\lvert \bn\rvert$ different $j$ such that $i_j^k>0$. By the previous comment, this is equivalent to select $\lvert \bn \rvert$ distinct levels  $j\leq N$.  From them, we need to choose  the $\bn(1)$ levels with of type~$1$, that is $\mathsf{g}(\lbl{j}{t}{\bx})=1$;  $\bn(2)$ levels of type~$2$; and so on. 
Let $q_l^k$  be the level of the $l$-th selected  type~$k$, $l\leq \bn(k)$, $k\in [d]$. 
By exchangeability, we rewrite~\eqref{eq:auxlem31} as
	\begin{align*}
		&\frac{\bn(1)!\cdots \bn(d)!}{N^{\lvert \bn \rvert}}\sum_{\substack{q_{1}^k,\ldots,q_{\bn(k)}^k\in[N]\\
		k\in[d],q_l^k\text{ all distinct}}} \E\bigg[ \prod_{k=1}^d \Big( \ind{\mathsf{g}(\lbl{q_{1}^k}{t}{\bx})=k,\dots,\mathsf{g}\left(\lbl{q_{\bn(r)}^k}{t}{\bx}\right)=k} \Big)\bigg]\\
		&=\frac{\bn(1)!\cdots \bn(d)!}{N^{\lvert \bn \rvert}} \binom{N}{\lvert \bn\rvert}\binom{\lvert \bn\rvert }{\bn(1),\ldots,\bn(d)}\P\left(E_{\bx,\bn,t}\right)\\
		&=\frac{N\cdots (N-\lvert \bn \rvert +1) }{N^{\lvert \bn\rvert}}
\P\left(E_{\bx,\bn,t}\right),
	\end{align*}
where $E_{\bx,\bn,t}$ is defined in \eqref{eq:Exnt}. The previous equation converges to the right-hand side of~\eqref{eq:prob} as $N\to\infty$. 
	
	Next, we show that the terms in~\eqref{eq:upperlower} corresponding to $\ell<\lvert \bn \rvert $ vanish as $N\to\infty$. To this end, fix some $\ell<\lvert \bn \rvert$ and consider the corresponding summand in~\eqref{eq:upperlower}.  As before, for a positive contribution it is necessary that if $i_j^k>0$ and $\mathsf{g}(\lbl{j}{t}{\bx})=k$, then $i_j^l=0$ for $l\neq k$. 

	In particular, we first need to select among~$N$ levels $\lvert \bn \rvert$ levels, with each level $j$ having $i_j^k>0$ for some~$k$. But in contrast to $\ell=\lvert \bn \rvert$, some levels can be chosen more than once and the number of \emph{distinct} lines is $\ell$.
For any $k\leq d$, let   $\bn^\star(k)=| \{i_j^k: i_j^k\neq 0, j\geq N\}|$ and  $q_l^k$ be the level of the $l$-th selected type~$k$, $l\leq \bn^\ast(k)$, $k\in [d]$. 
 Thus, we can upper bound the term in \eqref{eq:upperlower} corresponding to~$\ell$ by
	\begin{align*}
		& \frac{\bn(1)!\cdots \bn(d)!}{N^{\lvert \bn \rvert}}\sum_{\bn^\ast(1)=1}^{\bn(1)}\hspace{-.2cm}\dots \hspace{-.2cm} \sum_{\bn^\ast(k)=1}^{\bn(k)} 
\ind{ \ell=|\bn^\ast|}
\hspace{-.2cm}
\sum_{\substack{q_{1}^k,\dots,q_{\bn^\star(k)}^k\in[N]\\
k\in[d],q_l^k\text{ all distinct}}}
\hspace{-.2cm}
\E\bigg[  \prod_{k=1}^d \ind{\mathsf{g}(\lbl{q_{1}^k}{t}{\bx})=k,\dots,\mathsf{g}\left(\lbl{q_{\bn^\star(k)}^k}{t}{\bx}\right)=r}\bigg]\\
&\leq \frac{\bn(1)!\cdots \bn(d)!}{N^{\lvert \bn \rvert}}
\sum_{\bn^\ast(1)=1}^{\bn(1)}\dots \sum_{\bn^\ast(k)=1}^{\bn(k)}
\ind{ \ell=|\bn^\ast|}
 \binom{N}{\ell}\binom{\ell}{\bn^\star(1),\ldots,\bn^\star(d)}\\
&\leq \bn(1)!\cdots \bn(d)!\bn(1)\cdots \bn(d)\frac{N^{\ell}}{N^{\lvert \bn\rvert}},
	\end{align*}
 where we have use exchangeability and that the $\{q_l^k: l\leq \bn^\ast(k), k\leq d \}$ levels can be taken in the following form: first, we select $\ell$ distinct levels from the $N$ available and then, we choose the $\bn^\ast(k)$ levels of type~$k$, for every $k\leq d$.  Since $\ell <\lvert \bn\rvert$, the last equation converges to~$0$ as $N\rightarrow \infty$. 
\end{proof}

We are in position to prove the convergence theorem.

\begin{proof}[Proof of Theorem~\ref{thm:convergence}]
	Fix~$N\in\N$. By exchangeability, using the standard argument of applying a uniform permutation at every event in the $N$-lookdown, $\bX^{\bx,N}$ has the same distribution as the $\dtype$-type Moran model with large offspring.
	By~\cite[Ch. VI, Prop. 2.2]{Jacod1987}, for proving convergence of the $d$-type Moran model, it suffices to prove convergence for the $d$ one-dimensional marginals. 
    (To verify the condition in the proposition, we notice that all marginals jump due to the same Poisson processes and we can thus, in their notation, choose $t_n=t$.)
	Each of these marginals is a Moran model with mutation and large offspring that converges weakly to the $\Lambda$-Wright-Fisher process with mutation~\cite[Proposition 4.1.]{Etheridge2010}.
	In combination with \cite[Ch. VI, Prop. 2.2]{Jacod1987} it  thus follows that $\bX^{\bx,N}$ converges weakly to a process $\tilde{\bX}^{\bx}$.
	Next, we prove that for fixed $t$, the moments of $\tilde{\bX}_t^{\bx}$ are equal to the moments of $\bX_t^{\bx}$ (the process associated to~\eqref{eq:wfgenerator}).
	Fix $\bn\in \N_{0}^{\dtype}$. 
	Then, by dominating convergence $\E[(\tilde{\bX}_t^{\bx})^{\bn}]= \lim_{N\to\infty} \E[(\bX_t^{\bx,N})^{\bn}]$.
	By combining Lemmas~\ref{lem:asymmomentstypes} and \ref{lem:momentstypes}, we obtain $\lim_{N\to\infty}\E[(\bX_t^{\bx,N})^{\bn}]=\E\left[\bx^{\bA_t^{(\bn)}}\right].$
	Moreover, by the duality (Theorem~\ref{thm:duality}), $\E\left[\bx^{\bA_t^{(\bn)}}\right]=\E\left[(\bX_t^{\bx})^{\bn}\right]$.
	Altogether, this proves the result. 
\end{proof}

\section{Fixation line process and its explosion time}\label{sec:fixation_line_process}

To analyze the asymptotic behavior of the process $\bX^{\bx}$, 
the fixation lines associated with the lookdown construction are crucial.
To this end, we now prepare the definition and relevant results.

\begin{definition}[Fixation line]
For every $k\in \N_0$, the \emph{$k$-fixation line} is the process $F^k=(F^k_t:\, t\geq 0)$ taking values in $\N_0$ defined as 
	\[
F_t^k= \max\left\{n\in \N:\ \forall 0<m\leq n,\ \lbl{m}{t}{\bx} \text{ is measurable w.r.t. }\mathscr{F}_t\vee\sigma(\bU(j): 1\leq j\leq k) \right\},
\]
	with the convention that $\max\emptyset=0$ and 
\[F_t^0=\max\left\{n\in \N:\ \forall 0<m\leq n,\ \lbl{m}{t}{\bx} \text{ is measurable w.r.t. }\mathscr{F}_t\right\}.\]
	
\end{definition}
Informally, the $k$-fixation line at time~$t$ is the highest level such that the types of all individuals below it can be fully determined from the Poisson point processes (up to time~$t$) together with the values of the random variables $\bU(1), \ldots, \bU(k)$. 
For an illustration, see Figure~\ref{fig:fixation_line}. 

Without mutation ($\theta=0$), $F_t^0\equiv 0$ for all $t\geq 0$. 
In this context, Hénard \cite[Sect.~2.1]{henard2015} defined the $k$-th fixation line at time~$t$ as the smallest level~$\ell$ such that level $\ell+1$ carries (at time~$t$) an offspring of the individual occupying level~$k+1$ at time~$0$.
This definition coincides with ours. For an overview of the notion of fixation lines, see Hénard \cite[Sect.~1]{henard2015}, which also discusses related work by \cite{Labbe2014,Labbe2014a,Pfaffelhuber2006}.
The next proposition states the transition rates of the fixation lines, extending the results of Hénard  \cite[Lem.~2.3]{henard2015}.

\begin{proposition}[Fixation line rates]\label{prop:fixationlinerates}
	For each $k\in \N_0$, $F^k$ is a continuous-time Markov chain with $F^k_0=k$ that jumps from $n$ to $n+\ell$ at rate \[\ind{\ell=1}\bigg[\Lambda(\{0\})\binom{n+1}{2}+\theta (n+1)\bigg]+\binom{n+\ell}{\ell+1}\int_{(0,1]} r^{\ell+1}(1-r)^{n}\frac{\Lambda(\dd r)}{r^2},\qquad l\geq 1.\]
\end{proposition}

\begin{proof}
    Let $t\geq 0$, $k\in \N_0$. 
    First, we note that $F_t^k$ is piecewise constant and it may only change at the jumps of {\it{some}} Poisson point processes. More precisely, assume that $F_t^k=n$ for some $n\in \N$. Then, the first jump of $F^k$ after time $t$ is the first jump  after time $t$ of one of the following processes 
    \[
\{\mathcal{P}_{l,m}: l\leq n, m\leq n+1\}\cup\{\mathcal{N}\}\cup\{\mathcal{M}_l; l\leq n+1\}.
\]
    If the first jump comes from $\{\mathcal{P}_{l,m}; l\leq n, m\leq n+1\},$ then $F^k$ jumps from $n$ to $n+1$. This occurs at rate $\Lambda(\{0\})\binom{n+1}{2}$. 
    Now, suppose that  the first jump comes from $\mathcal{N}$, with an atom $(s,r,\bu)$. 
    If $\lvert\{m\leq n+1: \bu(m)\leq r \}\rvert \leq 1$, $F^k$ does not change. 
    Otherwise, $F^k$ jumps from $n$ to $n+\ell$ when $\bu(m)\leq r$ for exactly $\ell+1$ integers~$m$ with $m\leq n+\ell$ and $\bu(n+\ell+1)>r$. 
    The probability of such~$\bu$ is $r^{\ell+1}(1-r)^{n}$, because there are $n$ marked and $\ell+1$ unmarked positions in the first $n+\ell+1$ levels. There are $\binom{n+\ell}{\ell+1}$ ways of distribute the marks. 
    Integrating with respect to $r^{-2}\Lambda(\dd r)$ gives the rate. 
    Finally, if the first jump comes from $\{\mathcal{M}_l:\, l\leq n+1\}$, 
    then $F^k$ jumps from $n$ to $n+1$. This occurs at rate~$(n+1)\theta$.  
    Since these are a finite number of possibilities, the process $F_t^k$ is piecewise constant with exponential interarrival times, 
    we can conclude that it is a continuous-time Markov chain.
\end{proof}

Observe that for each $k$, the Markov chain $F^k$ stays at some state a positive amount of time, 
but it is possible that these holding times become infinitesimal and leading to infinitely many jumps within a finite time. This phenomenon is known as explosion. 

\begin{definition}[Explosion time]
	For each $k,n\in \N_0$ with $n\geq k$, the generalized right inverse of $F^k$ at~$n$ is defined as $T_n(k)\coloneqq \inf\{t\geq 0: F_t^k\geq n\}$. 
	Moreover, the explosion time of $F^k$ is \[T_{\infty}(k)\coloneqq \inf\{t\geq 0:\, F_t^k=\infty\},\]
	where $\inf \emptyset =\infty$.
\end{definition}
The condition of coming down from infinity, Equation \eqref{cond:cdi},   implies the almost sure explosion of the fixation lines. 
More precisely, for $k\in \N$,
\[
\{T_{\infty}(k)<\infty\}\subset\{\Lambda\text{-coalescent comes down from infinity}\}.
\]
This implication justifies the appearance of the coming down from infinity condition as a hypothesis in some of our results.

In the following lemma, we compute the mean explosion time in two particular cases.
\begin{proposition}\label{prop:explosionbeta}
Let $k\in \N$.
\begin{enumerate}
\item[(a)] \textbf{Wright-Fisher diffusion.} If $\theta\geq 0$ and $\Lambda=\Lambda(\{0\})\delta_0$, then
\[
    \E\left[T_{\infty}(k)\right]=2\sum_{j=k+1}^{\infty} \frac{1}{\Lambda(\{0\})j(j-1)+2j\theta }.
\]
\item[(b)] \textbf{Beta-coalescent.} If $\theta=0$ and $\Lambda$  given by~\eqref{eq:measurelambdabeta} with $\alpha\in (1,2)$, then \[\E\left[T_{\infty}(k)\right]=
		\alpha(\alpha-1)\int_{(0,1)} \frac{y^{k}(1-y)^{-1}}{(1-y)^{1-\alpha}-1}\dd y. 
	\]
\end{enumerate}
\end{proposition}

\begin{proof}
For $\theta\geq0$ and $\Lambda=\Lambda(\{0\})\delta_0$, then 
 $T_{\infty}(k)$ is the explosion time of a continuous-time Markov chain started at $k$ that jumps from $n$ to $n+1$ at rate  $\Lambda(\{0\})\binom{n+1}{2}+\theta (n+1)$. Then
$T_{\infty}(k)=\sum_{r=k}^{\infty}E_r,$ 
where $E_r\sim \mathrm{exp}\left(\Lambda(\{0\})\binom{r+1}{2} + \theta (r+1)\right)$ are independent exponential random variables. In particular,
\begin{equation*}
\E\left[T_{\infty}(k)\right]=2\sum_{j=k+1}^{\infty} \frac{1}{\Lambda(\{0\})j(j-1)+2j\theta }.
\end{equation*} 
Now, suppose that $\theta= 0$ and $\Lambda$ corresponds with a $\mathrm{Beta}(2-\alpha,\alpha)$-coalescent. By definition of the explosion time,
\[
	\E\left[T_{\infty}(k)\right]
        =\sum_{j= k+1}^{\infty} \left( \E\left[ T_j(k) \right] - \E\left[ T_{j-1}(k) \right] \right)
        =\sum_{j= k+1}^{\infty}\frac{\P(j-1\in\mathcal{S}_k)}{\Lambda_j},
\]
where $\mathcal{S}_k\coloneqq \left\{F^k_t:t\geq 0\right\}$ is the range of the fixation line started at~$k$, and $\Lambda_{j}$ is the rate at which \textit{a} fixation line jumps from $j-1$ to a higher level. Hénard \cite{henard2015} showed  that for $j\in \N$,
\begin{equation*}
	\Lambda_{j}=\frac{1}{\alpha\, \mathrm{Beta}(j-1,\alpha)} =\left(\alpha\int_{(0,1)}y^{j-2}(1-y)^{\alpha-1}\dd y\right)^{-1}.
\end{equation*}  
and that  the law of $\left\{F^k(t)-k:t\geq 0\right\}$ is independent of~$k$, see \cite[Eq.~(2.14) and sentence before Prop.~2.6]{henard2015}. We write $\mathcal{S}$ for a generic random set with this law.  Then, 
\begin{align*}
\E\left[T_{\infty}(k)\right]&=\sum_{j= k+1}^{\infty} \P(j-1-k\in\mathcal{S})\, \alpha\int_{(0,1)}y^{j-2}(1-y)^{\alpha-1}\dd y\\
&=\alpha\int_{(0,1)}y^{k-1}(1-y)^{\alpha-1}\sum_{j=0}^{\infty} \P(j\in\mathcal{S})y^j\dd y.
\end{align*}
Finally, by using that  $\sum_{j=0}^{\infty} \P(j\in\mathcal{S})y^j=(\alpha-1)y((1-y)-(1-y)^\alpha)^{-1}$, see \cite[Prop.~2.6]{henard2015}, we obtain the result.
        
        \end{proof}

\section{Paths coalescence time and proof of Theorem~\ref{thm:fixandstattimes}}\label{sec:proof_fixation_and_stationary_times}
To study how the progeny of a single individual contributes to the distribution of types,
we define for $t\geq 0$ and $l\in \N$, 
$\mathfrak{f}_t(l)$ as the proportion (or relative frequency) of the population at time $t$ that descendants from the individual on level~$l$ at time~$0$, i.e.
	\[\mathfrak{f}_t(l)\coloneqq \lim_{N\to\infty} \frac{1}{N} \sum_{m=1}^N \ind{\lbl{m}{t}{\bx}=(\bU(l),\bx)}, \qquad t\geq 0.
    \] 
    To avoid overloading the notation, and motivated by Theorem~\ref{thm:convergence}, we will abuse notation and write $\bX^{\bx}$ for $\lim_{N\to\infty}\bX^{\bx,N}$ throughout the remainder of the manuscript.
    Observe that then
\begin{equation}\label{eq:xtermp}
	\bX^{\bx}_t(i)=\sum_{l=1}^\infty \mathfrak{f}_t(l) \ind{\mathsf{g}( \bU(l),\bx)=i},\quad i\in[\dtype+1].
\end{equation} 

Since the random variables $(\bU(l): l\in \N)$ are independent and uniformly distributed in $[0,1]$, we have that almost surely for all $l\neq k$, $(\bU(l),\bx)\neq (\bU(k),\bx)$. 
It then follows from Condition~\eqref{cond:cdi}
that before the explosion time of the $k$-fixation line, 
the initial configuration $\lbl{l}{0}{\bx}=(\bU(l),\bx)$  contributes to $\bX^{\bx}$ for $j\leq k+1$; and after the explosion time, the initial configurations $\lbl{l}{0}{\bx}$ plays no role for any $l\geq k+1$.  
Note that the value $k+1$ is included in both cases. 
The next lemma makes this precise.

\begin{lemma}\label{lem:cruciallem}
Let $k\in \N$. 
	If Condition~\eqref{cond:cdi} holds, then almost surely  $\mathfrak{f}_t(l)>0$ for all $l\leq k+1$ and $t\in [0,T_{\infty}(k))$. Additionally, 
	$\mathfrak{f}_t(l)=0$ almost surely, for all $l\geq k+1$ and $t\in [T_{\infty}(k), \infty).$
\end{lemma}

\begin{proof}
    Fix $k\in \N$. We begin by proving the second statement. 
    Let $t\geq T_{\infty}(k)$. 
    By definition of $T_{\infty}(k)$, for all $m \in\N$,
    the function $\rho_t^{\bx}(m)$ is measurable with respect to  $\mathscr{F}_t\vee\sigma(\bU(l):  l \in  [ k ])$.
    Recall that the random variables $\{\bU(l): l\in \N\}$ are independent and uniformly distributed in $[0,1]$. Moreover, they are independent from the Poisson point processes $(\mathcal{M}_l(s): s\geq 0,  l\in \N)$.
    Then, almost surely for $m\in\N$
    \[
    \rho_t^{\bx}(m) \in \{ (\bU(l),\bx) : l \in [k] \} \cup  \{  (u,\bnu) : u\in [0,1]\}.
   \]
    Consequently, $\mathfrak{f}_t(l)=0$ for all $l\geq k+1$ so that the second claim of the lemma is proved. 
    
    For the first claim fix $t<T_{\infty}(k)$ 
    and let $l\leq k+1$.
    First, assume that $\theta=0$. 
    Embedded in the lookdown-construction, there is a coalescent process with values in the partition of~$\N$.
    We first describe the coalescent and then explain how it can be used for the proof.
    The coalescent starts at time~$t$ with infinitely many blocks $(B_m(t):m\in \N)$,
    each consisting of one singletons with $B_m(t)=\{m\}$,
    and runs backwards in time. 
    Block $B_m(t)$ sits on level $m$ at time~$t$.
    At time~$s\in[0,t]$, 
    if $\mathcal{P}_{l,m}$ has a jump at time~$s$, 
    the two blocks $B_l(s)$ and $B_m(s)$ coalesce.  
    In this way, if the jump is due to $\mathcal{P}_{l,m}$, 
    then $B_{k}(s-)=B_{k+1}(s)$ for all $k\geq m$, $B_l(s-)=B_l(s)\cup B_m(s)$, and the other blocks remain unchanged. 
    Similarly, if the jump is due to $\mathcal{N}$ with an atom $(s,r,(\bu(1),\bu(2),\ldots))$, then
    (recall, from Section~\ref{seq:lookdownconstruction}, $k_0=\min\{l:\, \bu(l)\leq r\}$ and) set $B_m(s-)=B_m(s)$ for all $m< k_0$,  $B_{k_0} (s-)=\bigcup_{l: \bu(l)\leq r}B_l(s)$, and all other blocks are relabeled so that $B_m(s-)$ contains the smallest element not in $B_1(s-),\dots, B_{m-1}(s-)$.
    In particular, the resulting coalescent is a standard $\Lambda$-coalescent.
    Now that the coalescent it described, we explain how it can be used for the proof.
    Observe that for $s\in [0,t]$, $\mathfrak{f}_s(m)=\lim_{N\to\infty}\frac{1}{N}\sum_{l=1}^N\ind{B_l(s)\subseteq B_m(0)}$.
    In particular, $\mathfrak{f}_t(m)=\lim_{N\to\infty}\frac{1}{N}\sum_{l=1}^N\ind{l\in B_m(0)}$.
    Condition~\eqref{cond:cdi} implies that the  $\Lambda$-coalescent comes down from infinity, i.e. there are only a finite number of blocks, for every $s\in[0,t)$.
    Because $t<T_k(\infty)$,
    $B_l(0)$ is non-empty for all $l\leq k+1$.
    Moreover, it is well-known that if a $\Lambda$-coalescent comes down from infinity, 
    then its asymptotic frequencies are proper, that is, $\sum_{l=1}^{\infty} \mathfrak{f}_t(l)=1$ almost surely~\cite[Thm 8]{Pitman1999}.
    It follows from \cite[Prop.~2.8]{Bertoin2006a}, that $f_t(l)=0$ implies either $B_l(t)$ is a singleton or empty.
    By the same result, if a coalescent is proper, none of the blocks is a singleton.
    Thus, all non-empty blocks then have positive asymptotic frequencies.
    Because $t<T_{\infty}(k)$, also $\mathfrak{f}_{t}(l)>0$ for $l\leq k+1$.
    This finishes the proof for $\theta=0$.
    
    For $\theta>0$, we can use the theory of distinguished coalescents developed in~\cite{Foucart2011}. A distinguished coalescent is a coalescent process with a distinguished block. Such a coalescent is embedded in the lookdown-construction if we add a $0$ level, and we prescribe that the $l$th block coalesces at time~$s$ with the distinguished block if $\mathcal{M}_l$ has a jump at time~$s$. (All other transitions are as before.) By Foucart \cite[Rem.~3.1]{Foucart2011}, all blocks different from the distinguished block have proper asymptotic frequencies in our context under Condition~\eqref{cond:cdi}. Moreover, this distinguished coalescent comes down from infinity under the same conditions as the one without the distinguished block~\cite[Thm. 4.1]{Foucart2011}. Thus, also in the case $\theta>0$, $\mathfrak{f}_s(l)>0$ for any $l\leq k+1$ and $s\in(0,t]$.
\end{proof}

With the previous lemma we can provide the coalesce time of two different initial configurations. More precisely, given two initial distributions of type frequencies $\bx,\by\in \Delta_{\dtype}$, we define the first level for which the types differ at time~$0$, i.e.
\[
D_{\bx,\by}\coloneqq \min\{i: \type{\bU(i)}{\bx}\neq\type{\bU(i)}{\by}\},
\] 
with the convention that $\min\emptyset = \infty$. Then we have the following result.

\begin{theorem}[Coalescence times] \label{thm:coaltimes}
	If Condition~\eqref{cond:cdi} holds, then almost surely, for all $\bx,\by\in \Delta_{\dtype}$,  the path coalescence time of the processes with paths starting at $\bx$ and $\by$ is given by
	\[
T_{\mathrm{coal}}^{\bx,\by}\coloneqq\inf\left\{t\geq 0: \bX_t^{\bx}=\bX_t^{\by}\right\}\overset{(d)}{=}T_{\infty}(D_{\bx,\by}-1).\]
\end{theorem}
\begin{proof}
	The main idea is that the type of the individual at  level~$D_{\bx,\by}$ at time~$0$ differs under~$\bx$ and~$\by$, 
	whereas individuals on the first $D_{\bx,\by}-1$ levels have the same type under~$\bx$ and~$\by$. 
	We will distinguish between $t\geq T_{\infty}(D_{\bx,\by}-1)$ and $t< T_{\infty}(D_{\bx,\by}-1)$. 
	In the former case, by Lemma~\ref{lem:cruciallem}, $\mathfrak{f}_t(k)=0$ for all $k\geq D_{\bx,\by}$ and thus $\bX_t^{\bx}=\bX_t^{\by}$, thanks to \eqref{eq:xtermp}. 
	For the latter case, recall the definition of the type function \eqref{eq:typingfct} and assume $\mathsf{g}( \bU(D_{\bx,\by}),\bx)=v$ and $ \mathsf{g}(\bU(D_{\bx,\by}),\by)=w$ for $v,w\in [\dtype]$ with $v\neq w$. 
    Without loss of generality, assume $v< w$.
    We will show that 
    \[
    \bX_t^{\bx}(1)+\cdots+\bX_t^{\bx}(v)\neq \bX_t^{\by}(1)+\cdots+\bX_t^{\by}(v),
    \] 
    which implies $\bX_t^{\bx}\neq \bX_t^{\by}$.
    By the assumption at the beginning of the proof, 
    \[\sum_{i=1}^{v}\by(i) <\bU({D_{\bx,\by}})\leq \sum_{i=1}^{v}\bx(i).\]
    This implies that, if $\mathsf{g}(u,\bx)> v$ for some $u\in [0,1]$ then $\mathsf{g}(u,\by)> v$.
    Or equivalently
    if $\mathsf{g}(u,\by)\leq v$, then also $\mathsf{g}(u,\bx)\leq v$.
    Additionally, $\ind{\{\mathsf{g}(D_{\bx,\by},\bx)\leq v\}}\ind{\{\mathsf{g}(D_{\bx,\by},\by)> v\}}=1$. 
	Thus, 
	\begin{align*}
		\bX_t^{\bx}(1)+\ldots+\bX_t^{\bx}(v)&= \sum_{k=1}^{\infty}\ind{\{\mathsf{g}(\bU(k),\bx)\leq v\}}\mathfrak{f}_t(i) \\
		&= \sum_{k=1}^{\infty}\ind{\{\mathsf{g}(\bU(k),\bx)\leq v\}}(\ind{\{\mathsf{g}(\bU(k),\by)\leq v\}}+\ind{\{\mathsf{g}(\bU(k),\by)> v\}})\mathfrak{f}_t(i)\\
		&\geq \bX_t^{\by}(1)+\ldots+\bX_t^{\by}(v)+  \mathfrak{f}_t(D_{\bx,\by})>\bX_t^{\by}(1)+\ldots+\bX_t^{\by}(v),
	\end{align*}

    This completes the proof.

\end{proof}

\begin{remark}[A comb-like metric]
    The previous result yields a random metric on $\Delta_{\dtype}$. 
    To this end, given $\bU=(\bU(i):i\in \N)$ on $[0,1]^{\infty}$, define 
\[
d_{\bU}(\bx,\by)\coloneqq T_{\infty}(D_{\bx,\by}-1) \quad \text{for }\bx,\by\in \Delta_d.
\]
    This indeed defines a metric on $\Delta_{\dtype}$. Moreover, it defines an ultrametric because the next inequality follows  $d_{\bU}(\bx,\bz)\leq \max\{d_{\bU}(\bx,\bz),d_{\bU}(\bz,\by)\}$.
    
    In the two types case $\dtype=1$, there appears to be a connection with the comb metric~\cite[Prop. 3.1]{Lambert2019}.
    A \emph{comb} is a function $f:[0,\infty)\to[0,\infty)$ such that for any $\varepsilon>0$, $\{f\geq \varepsilon\}$ is a finite set.
    The corresponding \emph{comb metric} $d_f(x,y)\coloneqq \sup_{z\in (x\wedge y,x\vee y)}f(z)$ defines an ultrametric distance on $\{f=0\}$.
    The function given by $C(x)=\sum_{k\geq 1}T_{\infty}(k)\ind{\bU(k)}(x)$ is a comb
    (since, because of~\eqref{cond:cdi}, there are only finitely many $k$ such that $T_{\infty}(k)>\varepsilon$ a.s., see proof of Theorem~\ref{thm:coaltimes}).
    In particular, for $x,y\in[0,1]$,  $d_C(x,y)=d_{\bU}(x,y)$.
    Note that if $(T_{\infty}(k):\, k\in \N)$ were independent (which they are not!), $d_{\bU}(x,y)$ would be the Kingman comb defined in Lambert and Schertzer \cite[Prop. 3.1]{Lambert2019}. 
\end{remark}

In the following subsections, we analyze the asymptotic behavior of the process $\bX^\bx$. Specifically, we prove Theorem~\ref{thm:fixandstattimes}(i) in the next subsection and part (ii) in the one after.

\subsection{Fixation times and proof of Theorem~\ref{thm:fixandstattimes}-(i)}\label{sec: fixtime}
Throughout this subsection, we assume~$\theta=0$, i.e. mutations  are absent. We begin by presenting the consequences of Theorem~\ref{thm:coaltimes} for fixation times.  
The following corollary is an immediate consequence of the theorem and is therefore stated without proof.

\begin{corollary}[Fixation time] \label{coro:fixationetimes}
	Assume $\theta=0$. If Condition~\eqref{cond:cdi} holds, then for all $\bx\in \Delta_{\dtype}$ and $i\in [\dtype+1]$, \[\inf\{t\geq 0:\, \bX_t^{\bx}=\e{i}\}\overset{(d)}{=}
	T_{\infty}(D_{\bx,\e{i}}-1).\]
\end{corollary}

Recall from \eqref{eq:defVktyp} that $\mathfrak{m}^{\bx}_n$ is the first 
index of the random uniform sequence~$\bU$ such that $g(\bU(k),\bx)=n$.
Similarly, $V_j^{\bx}$ is the first index such that a $j$th type appears for the first time;
or equivalently, it is the
level at time~$0$ such that exactly $j-1$ distinct types are present below that level, and the individual at that level carries a type not among those $j-1$ types.
 
We now turn to the proof of Theorem~\ref{thm:fixandstattimes}(i). Before giving the formal argument, we outline the main idea using Fig.~\ref{fig:lookdown}. While the figure shows the case $\theta>0$, we are in this subsection in the situation $\theta=0$ -- so the black arrows and lines, should not appear.
Then, if we are, for example,
interested in the first time the type that is colored blue disappears from the population, 
we need to find the first level (at time~$0$) where this color appears (in the case of the Figure, it is level $3$). 
Next, find the explosion time of the fixation line starting from the level just below (so level $2$ in the figure).
At that time, all the lines colored blue have been \emph{pushed out} from the population by offspring of the individuals below (level $1$ and $2$). 

\begin{proof}[Proof of Theorem~\ref{thm:fixandstattimes}-(i)]
    Fix $n\in [\dtype+1]$ and $\bx\in \Delta_{\dtype}$; recall that $\theta=0$. 
    For part (a), 
    first consider the case $\bx(n)=0$.
    Then, on the one hand $T_{\mathrm{lost},n}^{\bx}=0$, on the other hand $\mathfrak{m}^{\bx}_n=\infty$ and thus also $T_{\infty}(\mathfrak{m}^{\bx}_n-1)=0$, proving the result in that case.
    Next, assume $\bx(n)>0$.
    Note that by Lemma~\ref{lem:cruciallem} for $k\in[\mathfrak{m}^{\bx}_n]$, we have $\mathfrak{f}_t(k)>0$ for any $t<T_{\infty}(\mathfrak{m}^{\bx}_n-1)$.
    Thus, $\bX_t^{\bx}(n)>0$ for $t<T_{\infty}(\mathfrak{m}^{\bx}_n-1)$. 
    Then, by definition of $T_{\mathrm{lost},n}^{\bx}$ we have $T_{\infty}(\mathfrak{m}^{\bx}_n-1)\leq T_{\mathrm{lost},n}^{\bx}$.
    On the other hand, by Lemma~\ref{lem:cruciallem} for any $t\geq T_{\infty}(\mathfrak{m}^{\bx}_n-1)$ and $k\geq \mathfrak{m}^{\bx}_n$, we have $\mathfrak{f}_t(k)=0$.
    Hence, $\bX_t(i)=0$ for all $i\neq n$ and $t\geq T_{\infty}(\mathfrak{m}^{\bx}_n-1)$.
    By definition of $T_{\mathrm{lost},n}^{\bx}$ we have $T_{\infty}(\mathfrak{m}^{\bx}_n-1) \geq T_{\mathrm{lost},n}^{\bx}$ and the result follows.
    
    For part (b), if $\lvert \{i:\bx(i)>0\}\rvert\leq n$, 
    then on the one hand $T_{\mathrm{fix},n}^{\bx}=0$, 
    on the other hand, $V_{n+1}^{\bx}=\infty$ and thus $T_{\infty}(\infty)=0$, 
    proving the result in that case.
    Thus we now assume $\lvert \{i:\bx(i)>0\}\rvert> n$. Set $J=\left\{j\in [\dtype+1]: \mathsf{g}(\bU(k),\bx)=j\text{ for some }k\in [V_{n+1}^{\bx}-1]\right\}.$ Note that $\lvert J\rvert = n$.
    By Lemma~\ref{lem:cruciallem} for any $k\in[V_{n+1}^{\bx}]$, we have $\mathfrak{f}_t(k)>0$ for any $t<T_{\infty}(V_{n+1}^{\bx}-1)$. 
    Thus, $\bX_t^{\bx}(j)>0$ for all $j\in J\cup \left\{\mathsf{g}(\bU( V_{n+1}^{\bx}),\bx)\right\}$ and $t<T_{\infty}(V_{n+1}^{\bx}-1)$. Then, by definition of $T_{\mathrm{fix},n}$ we have $T_{\infty}(V_{n+1}^{\bx}-1)\leq T_{\mathrm{fix},n}^{\bx}$.
    On the other hand, by Lemma~\ref{lem:cruciallem} for any $t\geq T_{\infty}(V_{n+1}^{\bx}-1)$ and  $k\geq V_{n+1}^{\bx}$, we have $\mathfrak{f}_t(k)=0$.
    Hence, $\bX_t^{\bx}(k)=0$ for all $k\in [\dtype]\setminus J$ and $t\geq T_{\infty}(V_{n+1}^{\bx}-1)$. By definition of $T_{\mathrm{fix},n}$ we have $T_{\infty}(V_{n+1}^{\bx}-1) \geq T_{\mathrm{fix},n}^{\bx}$.
    
    This completes the proof of part (i) of the theorem.
\end{proof}

To obtain the explicit formulas in Corollary~\ref{coro:fixandstattimes_explicit}, we require a handle on the distribution of~$V_n^{\bx}$.
\begin{lemma}[Coupon collector waiting time distribution] \label{lem:coupon}
	For all $\bx \in \mathrm{int}(\Delta_{\dtype})$, $n\in [\dtype]$ and $k\geq 2$, 
    \[
    \P(V_{n+1}^{\bx}=k)=\sum_{\ell=1}^n (-1)^{n-\ell}\binom{\dtype-\ell}{d-n} \sum_{1\leq c_1<\cdots<c_{\ell}\leq \dtype +1} \Bigg(\sum_{i=1}^\ell \bx(c_i)\Bigg)^{k-1}\Bigg(1-\sum_{i=1}^\ell \bx(c_i)\Bigg)
    \]
\end{lemma}

\begin{proof}
	Fix $\bx\in \mathrm{int}(\Delta_{\dtype})$, $n\in[\dtype]$, and $k\geq 2$.
	The event $\{V_{n+1}^{\bx}=k\}$ can be rephrased in terms of a coupon collector problem with non-uniform collection probabilities. 
	To be more specific, this relates to a scenario involving $\dtype+1$ coupons, which are numbered from $1$ to $\dtype+1$. 
	At each discrete time step, one coupon is received, and it is coupon~$i$ with probability $\bx(i)$.
	For $m\in [\dtype+1]$, let $q(s,m)$ be the probability that $m$th last coupon is received for the first time at time~$s$; 
	so $q(s,1)$ is the probability for the last coupon being received at time~$s$, 
	$q(s,2)$ is the probability for the second last coupon being received at time~$s$, and so on.
	It is well-known~\cite[Eq. (2)]{Schelling1954} that
	\begin{equation}
		q(s,m)=\sum_{l=m-1}^{\dtype-1} (-1)^{l-m+1}\binom{l}{m-1}\sum_{1\leq c_1<\ldots<c_{l+1}\leq \dtype+1} \Big(\sum_{j=1}^{l+1}\bx(c_j)\Big)\Big(1-\sum_{j=1}^{l+1}\bx(c_j)\Big)^{s-1}. \label{eq:schelling}
	\end{equation}
	The time the $(d+1-n)$th last coupon is received (so $m = d+1-n$) is also the time that the $(n+1)$th first is received (for $n=d$ this means that the time the last coupon is received is also the first time the $(d+1)$th coupon is received).
	In particular, $q(s,d+1-n)=\P(V_{n+1}^{\bx}=s)$.
	Thus, using~\eqref{eq:schelling} with $m=d+1-n$ and summing in the second sum over the types that are collected in the first $s-1$ levels instead of the types that are not collected, we obtain
	\begin{equation}
		\sum_{l=d-n}^{d-1} (-1)^{l-d+n} \binom{l}{d-n}\sum_{1\leq c_1<\ldots<c_{d-l}\leq \dtype+1} \Big(1-\sum_{i=1}^{d-l}\bx(c_i)\Big)\Big(\sum_{i=1}^{d-l}\bx(c_i)\Big)^{s-1}
	\end{equation}
	Finally, by setting $\ell=d-l$ yields the result.
\end{proof}

We now proceed to the proof of Corollary~\ref{coro:fixandstattimes_explicit}.

\begin{proof}[Proof of Corollary~\ref{coro:fixandstattimes_explicit}]
Fix $\bx\in \mathrm{int}(\Delta_{\dtype})$ and $n\in[\dtype]$.
By Theorem~\ref{thm:fixandstattimes}-i(b) and the independence of $(\bU(k): k\in \N)$ and $\mathbb{F}:=(\mathscr{F}_t, t\geq 0)$, we have
\begin{equation}\label{eq:independence}
\E[T_{\mathrm{fix},n}]=\sum_{k=2}^{\infty} \P(V_{n+1}^{\bx}=k)\E[T_{\infty}(k-1)], \qquad n\in[\dtype+1].
\end{equation}
We first prove part (a) and suppose that $\Lambda_0=\Lambda(\{0\})\delta_0$.
By Proposition~\ref{prop:explosionbeta} we have that 
$\E\left[T_{\infty}(k)\right]=2(\Lambda(\{0\})k)^{-1}$. Then, by  Lemma~\ref{lem:coupon}, we obtain  
	\begin{align*}
		\E&[T_{\mathrm{fix},n}]\\
&=\frac{2}{\Lambda(\{0\})}\sum_{k=2}^{\infty} 
\sum_{\ell=1}^n (-1)^{n-\ell}\binom{\dtype-\ell}{d-n}
\hspace{-.2cm} \sum_{1\leq c_1<\cdots<c_{\ell}\leq \dtype +1}\hspace{-.3cm} \frac{\Big(1-\sum_{j=1}^\ell \bx(c_j)\Big)\Big(\sum_{j=1}^\ell \bx(c_j)\Big)^{k-1}}{k-1}\\
		&= -\frac{2}{\Lambda(\{0\}) }\sum_{\ell=1}^n (-1)^{n-\ell} \binom{d-l}{d-n} \sum_{1\leq c_1<\ldots <c_\ell\leq \dtype+1} \hspace{-.2cm}\Big(1-\sum_{j=1}^\ell\bx(c_j)\Big)\log\Big(1-\sum_{j=1}^\ell \bx(c_j)\Big),
	\end{align*}
	where we use that $-\log(1-a)=\sum_{k\geq 1}a^k/k$.  
    
    For the characteristic function of $T_{\mathrm{fix,n}}^{\bx}$, we use analogous arguments. Recall that $T_{\infty}(k)=\sum_{r=k}^{\infty}E_r,$
where $E_r\sim \mathrm{exp}\left(\Lambda(\{0\})\binom{r+1}{2} \right)$ are independent exponential random variables. Then, by Lemma \ref{lem:coupon} and Fubini's Theorem
\begin{align*}
	\E&[e^{itT_{\mathrm{fix,n}}^{\bx}}]=\sum_{k=2}^{\infty} \P(V_{n+1}^\bx=k)\E[e^{it T_{\infty}(k-1)} ]
=\sum_{k=1}^{\infty} \P(V_{n+1}^{\bx}=k+1) \prod_{r=k}^{\infty} \E[e^{itE_r}]\\
&= \sum_{\ell=1}^{n} (-1)^{n-\ell} \binom{d-\ell}{d-n} \\
    &\hspace{.5cm} \cdot \sum_{1\leq c_1<\ldots<c_\ell\leq d+1} \Bigg(1-{ \sum_{j=1}^{\ell}}\bx(c_j)\Bigg)\sum_{k=1}^{\infty}\Bigg(\sum_{j=1}^{\ell}\bx(c_j)\Bigg)^k\prod_{r=k}^\infty \Big(1-\frac{2it}{\Lambda(\{0\})(r+1)r} \Big)^{-1}.	
\end{align*}

Next, we prove part (b) and suppose that $\Lambda$ has Beta$(2-\alpha,\alpha)$ distribution.
We combine Equation~\eqref{eq:independence}, Proposition~\ref{prop:explosionbeta} and Lemma~\ref{lem:coupon} to obtain
	\begin{align*}
		\frac{\E[T_{\mathrm{fix},n}]}{\alpha(\alpha-1)}=& \sum_{\ell=1}^n (-1)^{n-\ell} \binom{d-\ell}{d-n}\\
&\cdot\sum_{1\leq c_1<\ldots <c_\ell\leq \dtype+1} \Big(1-\sum_{j=1}^\ell \bx(c_j)\Big) 
        \int_{(0,1)} \sum_{k= 1}^{\infty} \Big(y\sum_{j=1}^\ell \bx(c_j)\Big)^k \frac{ (1-y)^{-1}}{(1-y)^{1-\alpha}-1} \dd y.
	\end{align*}
By simplifying the geometric series, we obtain the result. 
\end{proof}

We now turn to the proof of the order in which the types disappear as stated in Proposition~\ref{prop:disappearanceorder}. 

\begin{proof}[Proof of Proposition~\ref{prop:disappearanceorder}]
Fix $\bx\in \Delta_{\dtype}$.
	We claim that \begin{equation}
		T_{\mathrm{lost},{c_{\dtype+1}}}^{\bx}<\ldots<T_{\mathrm{lost},{c_{1}}}^{\bx}\quad\text{ if and only if }\quad \mathfrak{m}^{\bx}_{c_{1}}<\ldots<\mathfrak{m}^{\bx}_{c_{\dtype+1}}, \label{eq:Tandm}
	\end{equation}
where $\mathfrak{m}^{\bx}_i$ is given by~\eqref{eq:defmx}.
Indeed, assume that $\mathfrak{m}^{\bx}_{i}>\mathfrak{m}^{\bx}_{j}$ for some $i\neq j$.
	By Lemma~\ref{lem:cruciallem}, for any $t\in[T_{\infty}(\mathfrak{m}^{\bx}_{i}-1),T_{\infty}(\mathfrak{m}^{\bx}_{j}-1))$, 
	we have $\mathfrak{f}_t(\mathfrak{m}^{\bx}_{j})>0$, 
	and $\mathfrak{f}_t(\ell)=0$ for all $\ell\geq \mathfrak{m}^{\bx}_{i}$. 
	Thus, by equation \eqref{eq:xtermp}
    \[\bX_t^\bx(i)=\sum_{\ell=1}^{\mathfrak{m}^{\bx}_{i}-1} \mathfrak{f}_t(\ell)\ind{\mathsf{g}(\bU(\ell),\bx)=i}=0 \qquad \mbox{ and } \qquad  \bX_t^\bx(j)>\mathfrak{f}_t(\mathfrak{m}^{\bx}_{j})>0.
    \]
Thus, $T_{\mathrm{lost},{k}}^{\bx}<T_{\mathrm{lost},{l}}^{\bx}$.
	For the other direction, note that $\mathfrak{m}^{\bx}_{i}\leq \mathfrak{m}^{\bx}_{j}$ implies $\mathfrak{m}^{\bx}_{j}<\mathfrak{m}^{\bx}_{i}$, 
	and thus by the previous computation $T_{\mathrm{lost},{i}}^{\bx}>T_{\mathrm{lost},{j}}^{\bx}$,
	which then implies $T_{\mathrm{lost},{i}}^{\bx}\geq T_{\mathrm{lost},{j}}^{\bx}$. Then~\eqref{eq:Tandm} is true.
    
In order to prove our result, we use the independence of $(\bU(k):\, k\in \N)$ and the geometric series, 
	\begin{align*}
		\P(T_{\mathrm{lost},{c_{\dtype+1}}}^{\bx}<\ldots<T_{\mathrm{lost},{c_{1}}}^{\bx})&=\P(\mathfrak{m}^{\bx}_{c_{1}}<\ldots<\mathfrak{m}^{\bx}_{c_{\dtype+1}})\\
		&= \sum_{n_2,\ldots,n_{\dtype+1}\geq0 } \bx(c_1)\prod_{i=2}^{\dtype+1}\Big\{ \Big(\sum_{j=1}^{i-1}\bx(c_j) \Big)^{n_i}\bx(c_i)\Big\}\\
		&=\bx(c_1)\prod_{i=2}^{\dtype+1} \frac{\bx(c_i)}{1-\sum_{j=1}^{i-1}\bx(c_j)}.
	\end{align*}
\end{proof}

Using the previous theorem and the following lemma, we will be able to prove our second main result, Theorem \ref{thm:meanfixk}. We begin by stating the lemma and providing its proof. It describes the distribution of $V_n^{\bx}$, which can be interpreted as the waiting time in a discrete-time coupon collector problem with non-uniform collection probabilities.

\subsection{Stationary times and proof of Theorem~\ref{thm:fixandstattimes}-(ii)}\label{sec:proofmixingtimes}
In this section, we prove that the explosion time of the $0$-fixation line is a stationary time for the process~$\bX$.
To this end, we rely on the following two lemmas. Throughout this subsection, we assume $\theta > 0$.

\begin{lemma}\label{lem:exploind}
	For all $\bx\in \Delta_{\dtype}$, $X_{T_{\infty}(0)}^{\bx}$ is independent from~$\bx$.
\end{lemma}
\begin{proof}
	By Lemma \ref{lem:cruciallem}, $\mathfrak{f}_t(k)=0$ for every $k\geq 1$ and $t>T_{\infty}(0)$. Then,  $\bX_{t}^{\bx}$ is independent of the value~$\bx$ for every 
		$t>T_{\infty}(0)$. Finally, we use the right continuity of the process to get the result.
\end{proof}
Because of Lemma~\ref{lem:exploind}, we can  write $\pi(\dd \by)\coloneqq\P\left(\bX_{T_{\infty}(0)}^{\bx}\in \dd \by\right)$ for an arbitrary $\bx\in \Delta_{\dtype}$.
\begin{lemma}\label{lem:invdist}
For any $\bx$, the distribution $\pi$ is stationary for~$\bX^{\bx}$.
\end{lemma}
\begin{proof}
	Let $\bx\in \Delta_{\dtype}$ arbitrary. Using the definition of $\pi$, the Chapman-Kolmogorov property,  and Lemma~\ref{lem:exploind}, we compute
	 \begin{align*}
		\P\left(\bX_t^\pi\in \dd \bz\right)&=\int_{\by\in\Delta_d}\P\left(\bX_t^\by\in \dd \bz\right)\pi(\dd \by)=\int_{\by\in\Delta_d}\P\left(\bX_t^\by\in \dd \bz\right)\P\left(\bX_{T_{\infty}(0)}^{\bx}\in \dd \by \right)\\
	 &=\P\left(\bX_{t+T_{\infty}(0)}^{\bx}\in \dd \bz \right)= \int_{\by\in\Delta_d}\P\left(\bX_{T_{\infty}(0)}^{\by} \in \dd \bz\right)\P\left(\bX_t^\bx\in \dd \by \right) \\
	 &= \pi(\dd \bz) \int_{\by\in\Delta_d}\P\left(\bX_t^\bx\in \dd \by \right)=\pi(\dd \bz).
	\end{align*}
\end{proof}

\begin{lemma}\label{lem:wfstrong}
	If $\Lambda=\Lambda(\{0\})\delta_0$, then $\bX_{T_{\infty}(0)}^{\bx}$ is independent of $T_{\infty}(0)$.
\end{lemma}

\begin{proof}
Recall that 
$T_{\infty}(0)=\sum_{k=0}^{\infty}E_k,$
where $E_k\sim \mathrm{exp}\left(\Lambda(\{0\})\binom{k+1}{2}+(k+1)\theta \right)$ are independent random variables  and  $E_k$ represents the occupation time of $F^0$ at  level~$k$.  Define for $k\in \N_0$,
	\[Y_k=\begin{cases}
		(l,m)\in \N^2, & \text{if }E_k\text{ increases because of an arrival of }\mathcal{P}_{l,m},\\
		(l,u)\in \N\times [0,1], & \text{if }E_k\text{ increases because of an arrival }\text{ of }\mathcal{M}_l\text{, say } (t,u).
	\end{cases}\]
	Then $Y_0,Y_1,\ldots$ and $E_0,E_1,\ldots$ are all mutually independent (one way to see this is via the colouring theorem of Poisson processes).
    Moreover, $\bX_{T_{\infty}(0)}^{\bx}$ is a measurable function of $(Y_k:\, k\in \N_0)$.
    In particular, it is independent of $E_0,E_1,\ldots$ and thus of $T_{\infty}(0)$.
\end{proof}
\begin{remark}
	If $\Lambda((0,1])>0$, then $F^0$ does not necessarily visit every state. Thus, it is plausible that the explosion time is not independent from the type distribution at the time of explosion. However, we did not prove this.
\end{remark}

We are ready to prove the remaining part  of Theorem~\ref{thm:fixandstattimes}.
In this case, we will show that $T=T_{\infty}(0)$ is stationary time with $\pi$ is stationary distribution.
If in addition  $\Lambda=\Lambda(\{0\})\delta_0$, then $T$ is a strong stationary time.

\begin{proof}[Proof of Theorem~\ref{thm:fixandstattimes}-(ii)]
	Recall that $\pi$ is an invariant distribution by Lemma~\ref{lem:invdist}.
    Using  Condition~\eqref{cond:cdi}, we know that $T_{\infty}(0)$ is almost surely finite.
    Moreover, $T_{\infty}(0)$ is a stopping time for the filtration $\mathscr{F}_t$ given at \eqref{eq:Poifiltration}.
    Therefore, $T_{\infty}(0)$ is a stationary time and by  Proposition \ref{prop:fixationlinerates} we obtain the first part.
    Now suppose that $\Lambda=\Lambda(\{0\})\delta_0$.
    Lemma~\ref{lem:wfstrong} then implies that  it is a strong time.
    
\end{proof}

Finally, we can prove the mean strong stationary time formula in the case of the Wright-Fisher diffusion.
\begin{proof}[Proof of Corollary~\ref{coro:meanstrongstationarytime}]
    By  Theorem~\ref{thm:fixandstattimes}-ii, $T_{\infty}(0)$ is a strong stationary time under the conditions of the corollary.
    By Proposition~\ref{prop:explosionbeta},
	\begin{align*}
    \E\left[T_{\infty}(0)\right]=2\sum_{j=1}^{\infty} \frac{1}{\Lambda(\{0\})j(j-1)+2j\theta }.
    \end{align*} 
    The result holds by using that $\sum_{i=1}^{\infty} \frac{a}{i(i+a)}=\psi(1+a)+\gamma$, see Abramowitz and Stegun \cite[Eq. 6.3.16]{Abramowitz1972}.  
\end{proof}

\section*{Acknowledgement}
Sebastian Hummel was funded by the Deutsche Forschungsgemeinschaft (DFG, German Research Foundation) -- Projektnummer 449823447. The authors thank the Hausdorff Research Institute for Mathematics in Bonn, where part of this research was carried out during the 2022 Junior Trimester Program ``Stochastic modelling in the life science: From evolution to medicine".

\bibliographystyle{amsplain}

\begin{thebibliography}{99}

\bibitem{Abramowitz1972}
Milton Abramowitz and Irene~A Stegun, \emph{Handbook of mathematical functions
  with formulas, graphs, and mathematical tables. national bureau of standards
  applied mathematics series 55. tenth printing.}, ERIC, 1972.

\bibitem{Berestycki2009}
Nathana\"{e}l Berestycki, \emph{Recent progress in coalescent theory}, Ensaios
  Matem\'{a}ticos, vol.~16, Sociedade Brasileira de Matem\'{a}tica, Rio de
  Janeiro, 2009.

\bibitem{Bertoin2006a}
Jean Bertoin, \emph{Random fragmentation and coagulation processes}, vol. 102,
  Cambridge University Press, 2006.

\bibitem{Birkner2009}
Matthias Birkner, Jochen Blath, Martin M{\"o}hle, Matthias Steinr{\"u}cken, and
  Johanna Tams, \emph{A modified lookdown construction for the
  {Xi}-{Fleming}-{Viot} process with mutation and populations with recurrent
  bottlenecks}, ALEA Lat. Am. J. Probab. Math. Stat. \textbf{6} (2009), 25--61
  (English).

\bibitem{Cordero2022a}
Fernando Cordero, Adri{\'a}n~Gonz{\'a}lez Casanova, Jason Schweinsberg, and
  Maite Wilke-Berenguer, \emph{$\lambda$-{coalescents arising in a population
  with dormancy}}, Electron. J. Probab. \textbf{27} (2022), no.~none, 1 -- 34.

\bibitem{Cordero2022}
Fernando Cordero, Sebastian Hummel, and Emmanuel Schertzer, \emph{{General
  selection models: Bernstein duality and minimal ancestral structures}}, Ann.
  Appl. Probab. \textbf{32} (2022), no.~3, 1499 -- 1556.

\bibitem{Diaconis1990}
Persi Diaconis and James~Allen Fill, \emph{Strong stationary times via a new
  form of duality}, Ann. Probab. \textbf{18} (1990), no.~4, 1483--1522.
  \MR{1071805}

\bibitem{Donnelly1999}
P.~Donnelly and T.~G. Kurtz, \emph{Particle representations for measure-valued
  population models}, Ann. Probab. \textbf{27} (1999), 166--205.

\bibitem{donnelly2000convergence}
Peter Donnelly and Eliane~R Rodrigues, \emph{Convergence to stationarity in the
  {M}oran model}, J. Appl. Probab. \textbf{37} (2000), no.~3, 705--717.

\bibitem{Durrett2008}
Richard Durrett, \emph{Probability models for {DNA} sequence evolution}, second
  ed., Probability and its Applications (New York), Springer, New York, 2008.
  \MR{2439767}

\bibitem{Eldon2006}
B.~Eldon and J.~Wakeley, \emph{Coalescent processes when the distribution of
  offspring number among individuals is highly skewed}, Genetics \textbf{172}
  (2006), 2621--2633.

\bibitem{Etheridge2010}
A.~M. Etheridge, R.~C. Griffiths, and J.~E. Taylor, \emph{A coalescent dual
  process in a {M}oran model with genic selection, and the lambda coalescent
  limit}, Theor. Popul. Biol. \textbf{78} (2010), 77--92.

\bibitem{Fill2016}
James~Allen Fill and Vince Lyzinski, \emph{Strong stationary duality for
  diffusion processes}, J. Theor. Probab. \textbf{29} (2016), no.~4,
  1298--1338. \MR{3571247}

\bibitem{Foucart2011}
Cl{\'e}ment Foucart, \emph{Distinguished exchangeable coalescents and
  generalized {Fleming-Viot} processes with immigration}, Adv. Appl. Probab.
  \textbf{43} (2011), no.~2, 348--374.

\bibitem{henard2015}
O.~Hénard, \emph{The fixation line in the ${\Lambda}$-coalescent}, Ann. Appl.
  Probab. \textbf{25} (2015), 3007--3032.

\bibitem{Jacod1987}
Jean Jacod and Albert~N. Shiryaev, \emph{Limit theorems, density processes and
  contiguity}, pp.~535--571, Springer Berlin Heidelberg, Berlin, Heidelberg,
  1987.

\bibitem{Labbe2014a}
Cyril Labb\'{e}, \emph{From flows of {$\Lambda$}-{F}leming-{V}iot processes to
  lookdown processes via flows of partitions}, Electron. J. Probab. \textbf{19}
  (2014), no. 55, 49. \MR{3227064}

\bibitem{Labbe2014}
\bysame, \emph{Genealogy of flows of continuous-state branching processes via
  flows of partitions and the {E}ve property}, Ann. Inst. Henri Poincar\'{e}
  Probab. Stat. \textbf{50} (2014), no.~3, 732--769. \MR{3224288}

\bibitem{Lambert2019}
Amaury Lambert and Emmanuel Schertzer, \emph{{Recovering the Brownian
  coalescent point process from the Kingman coalescent by conditional
  sampling}}, Bernoulli \textbf{25} (2019), no.~1, 148 -- 173.

\bibitem{liggett2010}
T.~M. Liggett, \emph{Continuous time markov processes: An introduction},
  American Mathematical Society, Providence, RI, 2010.

\bibitem{Littler1975}
R.~A. Littler, \emph{Loss of variability at one locus in a finite population},
  Math. Biosci. \textbf{25} (1975), no.~1-2, 151--163. \MR{392025}

\bibitem{Miclo2017}
Laurent Miclo, \emph{Strong stationary times for one-dimensional diffusions},
  Ann. Inst. Henri Poincar\'{e} Probab. Stat. \textbf{53} (2017), no.~2,
  957--996. \MR{3634282}

\bibitem{Oliver2012}
Matthew~K. Oliver and Stuart~B. Piertney, \emph{Selection maintains mhc
  diversity through a natural population bottleneck}, Molecular Biology and
  Evolution \textbf{29} (2012), no.~7, 1713--1720.

\bibitem{Pfaffelhuber2006}
P.~Pfaffelhuber and A.~Wakolbinger, \emph{The process of most recent common
  ancestors in an evolving coalescent}, Stochastic Process. Appl. \textbf{116}
  (2006), no.~12, 1836--1859.

\bibitem{Pitman1999}
J.~Pitman, \emph{Coalescents with multiple collisions}, Ann. Probab.
  \textbf{27} (1999), 1870--1902.

\bibitem{Sa99}
Serik Sagitov, \emph{The general coalescent with asynchronous mergers of
  ancestral lines}, J. Appl. Probab. \textbf{36} (1999), 1116--1125.
  \MR{1742154}

\bibitem{Schelling1954}
Hermann~Von Schelling, \emph{Coupon collecting for unequal probabilities},
  Amer. Math. Monthly \textbf{61} (1954), no.~5, 306--311.

\bibitem{schweinsberg2000}
J.~Schweinsberg, \emph{Coalescents with simultaneous multiple collisions},
  Electron. J. Probab. \textbf{5} (2000), no. 12, 50 pp.

\bibitem{sved2008divergence}
John~A Sved, Allan~F McRae, and Peter~M Visscher, \emph{Divergence between
  human populations estimated from linkage disequilibrium}, Am. J. Hum. Genet.
  \textbf{83} (2008), no.~6, 737--743.

\bibitem{Tavare1984}
S.~Tavar{\'e}, \emph{Line-of-descent and genealogical processes, and their
  applications in population genetics models}, Theor. Popul. Biol. \textbf{26}
  (1984), 119--164.

\end{thebibliography}

\end{document}